\documentclass[preprint,sort&compress,final,12pt]{elsarticle}

\usepackage{amsmath}
\usepackage{amssymb}
\usepackage{graphicx}
\usepackage{latexsym}
\usepackage{epstopdf}
\usepackage{natbib}

\newtheorem{theorem}{Theorem}[section]

\newtheorem{definition}[theorem]{Definition}
\newtheorem{lemma}[theorem]{Lemma}

\newtheorem{remark}[theorem]{Remark}
\numberwithin{equation}{section}

\def\Proof{\noindent{\bf Proof.}~}
\def\qed{\hfill$\square$\smallskip}
\def\dint{\displaystyle\int}
\def\dlim{\displaystyle\lim}

\def\dsum{\displaystyle\sum}
\def\dint{\displaystyle\int}

\def\dfrac#1#2{\frac{\displaystyle {#1}}{\displaystyle {#2}}}
\def\Re{\mathrm{Re}}

\def\Res{\mathrm{Res}}

\def\R{{\Bbb R}}
\def\C{{\Bbb C}}

\journal{Stochastic Processes and their Applications}

\begin{document}

\begin{frontmatter}

\title{Moment Boundedness of Linear Stochastic Delay Differential Equation with Distributed Delay}

\author[au1]{Zhen Wang}

\author[au1]{Xiong Li\footnote{Supported by the Natural Science Foundation of China (NSFC 11031002), the Fundamental
Research Funds for the Central Universities and the Scientific Research Foundation
for the Returned Overseas Chinese Scholars, State Education Ministry.}}

\address[au1]{School of Mathematical Sciences, Beijing Normal University, Beijing 100875, P.R. China.}

\author[au3]{Jinzhi Lei\footnote{Supported by the Natural Science Foundation of China  (NSFC 11272169).}}
\address[au3]{Zhou Pei-Yuan Center for Applied Mathematics, Tsinghua University, Beijing 100084, P.R. China.}

\begin{abstract}
This paper studies the moment boundedness of solutions of linear stochastic delay differential equations with distributed delay. For a linear stochastic delay differential equation, the first moment stability is known to be identical to that of the corresponding deterministic delay differential equation. However, boundedness of the second moment is complicated and depends on the stochastic terms. In this paper, the characteristic function of the equation is  obtained through techniques of Laplace transform. From the characteristic equation, sufficient conditions for the second moment to be bounded or unbounded are proposed.
\end{abstract}


\end{frontmatter}

\section{Introduction}
Time delays are known to be involved in many processes in biology, chemistry, physics, engineering, \textit{etc.},
and delay differential equations are widely used in describing these processes.  Delay differential equations have
been extensively developed in the past several decades (see \cite{Arino,Bellman,Hale93}). Furthermore, stochastic
perturbations are often introduced into these  deterministic systems in order to describe the effects of fluctuations
in real environment, and thus yield stochastic delay differential equations. Mathematically, stochastic delay
differential equations were first introduced by It\^{o} and Nisio in the 1960s \cite{Ito} in which the existence
and uniqueness of the solutions have been investigated. In the last several decades, numerous studies have been developed
toward the study of stochastic delay differential equations, such as stochastic stability, Lyapunov functional method,
Lyapunov exponent, stochastic flow, invariant measure, invariant manifold, numerical approximation and attraction \textit{etc.}
(see
\cite{Caraballo10,Prato96,Duan03,Duan04,Ivanov03,Khasminskii,Kuchler68,Kuchler00,Mandrekar94,Mao1,Mao2,Mao3,Mohammed92,Mohammed03,Mohammed04,Wu}
and the references therein). However, many basic issues remain unsolved even for a simple linear equation with constant coefficients.

In this paper, we study the following linear stochastic differential equation with distributed delay
\begin{equation}
\label{1.1}
\begin{array}{rcl}
dx(t)&=&\left(ax(t)+b\dint_{0}^{+\infty}K(s)x(t-s)ds\right)dt \\
&&{} + \left(\sigma_{0}+\sigma_{1}x(t)+\sigma_{2}\dint_{0}^{+\infty}K(s)x(t-s)ds\right)dW_t.
\end{array}
\end{equation}
Here $a,b$ and $\sigma_i(i=0,1,2)$ are constants, $W_t$ is a one dimensional Wiener process, and $K(s)$
represents the density function of the delay $s$. In this study, we always assume It\^{o} interpretation for the
stochastic integral. This paper studies the moment boundedness of the solutions of
\eqref{1.1}. Particularly, this paper gives the characteristic function of the equation, through
 which sufficient conditions for the second moment to be bounded or unbounded are obtained.

Despite the simplicity of \eqref{1.1}, which is a linear equation with constant coefficients, current understanding
for how the stability and moment boundedness depend on the equation coefficients is still incomplete. Most
of known results are obtained through the method of Lyapunov functional. The Lyapunov functional method
is useful for investigating the stability of differential equations, and has been well developed for
delay differential equations \cite{Hale93}, stochastic differential equations \cite{Mao1}, and stochastic
delay differential equations \cite{Khasminskii,Kuchler68,Mandrekar94,Mao1}. The Lyapunov functional method can usually give sufficient conditions for the stability of stochastic
delay differential equations. For general results one can refer to the Razumikhin-type theorems on the
exponential stability for the stochastic functional differential equation \cite[Chapter 5]{Mao1}.
However, these results often depend on the method of how the Lyapunov functional is constructed and are incomplete,
not always applicable for all parameter regions. For example, sufficient conditions for the $p^{\mathrm{th}}$ moment
stability of following stochastic differential delay equation
\begin{equation} \label{1.2}
dx(t)=\left(a x(t) + b x(t-\tau)\right)dt + \left(\sigma_{1}x(t) + \sigma_{2}x(t-\tau)\right)dW_t,
\end{equation}
can be obtained when $a<0$, but not for $a>0$ \cite[Example 6.9 in Chapter 5]{Mao1}.

In 2007, Lei and Mackey \cite{LeiM} introduced the method of Laplace transform to study the stability and
moment boundedness of the equation \eqref{1.1} with discrete delay ($K(s)=\delta(s-1)$). In this particular case, the characteristic equation was proposed, which yields a sufficient (and is also necessary
if not of the critical situation) condition for the boundedness of the second moment
(see Theorem 3.6 in \cite{LeiM}). This result gives a complete description for the second moment stability of the equation
\eqref{1.2} (the delay can be rescaled to $\tau$=1). Nevertheless, there is a disadvantage in the characteristic equation proposed in
\cite{LeiM} in that the characteristic function is not explicitly given by the equation coefficients. Therefore,
it is not convenient in applications.

The purposes of this paper are to study the stochastic delay differential equation \eqref{1.1} and to obtain a
characteristic function that is given explicitly through the equation coefficients.

Rest of this paper is organized as follows. In Section 2 we briefly introduce basic results for the fundamental
solutions of linear delay differential equations with distributed delay.
Main results and proofs of this paper are given in Section 3. In Section 3.1 we discuss the first moment stability and show that
it is identical to that of the unperturbed delay differential equation \eqref{2.1} (Theorem \ref{Thm3.3}). Section 3.2
focuses on the second moment. When the stochastic perturbation is an additive noise, the result is simple and the
bounded condition for the second moment is the same as the stability condition for the unperturbed delay differential
equation (Theorem \ref{Thm3.4}). However, in the presence of multiplicative white noises, boundedness for the second
moment depends on the perturbation terms. We prove that the second moment is unbounded provided that the trivial solution
of the unperturbed equation is unstable (Theorem \ref{Thm3.6}). When the trivial solution of the unperturbed equation is
stable, we obtain a characteristic equation, and the boundedness of the second moments depends on the maximum real parts
of all roots of the characteristic equation (Theorem \ref{Thm3.8}). The characteristic function is given explicitly through
the equation parameters. In Section 4, as applications, we obtain some useful conditions for the boundedness of the second
moment in some special situations (Theorem \ref{Thm4.1}), and also sufficient conditions for the second moment to be
unbounded (Theorem \ref{Thm4.2}). An example is studied in Section 5.

\section{Preliminaries}
In this section, we first give some basic results for the  fundamental solution
of a linear differential equation with distributed delay
\begin{equation}\label{2.1}
\dfrac{dx(t)}{dt}=ax(t)+b\int_{0}^{+\infty}K(s)x(t-s)ds.
\end{equation}
We also give sufficient and necessary conditions for the stability of the trivial solution of equation
(\ref{2.1}), which are useful for the rest of this paper. The linear delay differential equation has
been studied extensively and the existence and uniqueness of the solution can be referred to
\cite{Arino,Bellman,Hale93}.

First, we give some basic assumptions through out this paper. We always assume that the initial functions of \eqref{1.1} and \eqref{2.1} are $x=\phi\in BC\left((-\infty, 0], \R\right)$. Here $BC\left((-\infty, 0], \R\right)$ means the space of all bounded and continuous functions $\phi: (-\infty,0]\longrightarrow \R$
endowed with the norm
\begin{equation*}
\|\phi\|=\sup_{\theta\in(-\infty,0]}|\phi(\theta)|.
\end{equation*}
The delay kernel $K$
is a nonnegative piecewise continuous function defined on $[0, +\infty)$, satisfying
\begin{equation}\label{2.2}
\int_{0}^{+\infty} K(s)ds=1
\end{equation}
and there is a positive constant $\mu$ such that
\begin{equation}\label{2.3}
\int_{0}^{+\infty}e^{\mu s}K(s)ds< +\infty.
\end{equation}
We denote
\begin{equation}\label{2.4}
\rho =\int_{0}^{+\infty}e^{\mu s}K(s)ds
\end{equation}
for convenience.
For example, if we have gamma distribution delays:
\begin{equation} \label{2.5}
K(s)=\dfrac{r^{j}s^{j-1}e^{-r s}}{(j-1)!},~ s\geq 0,~ r>0,~ j=1,2,3,
\cdots
\end{equation}
then (\ref{2.2}) holds and (\ref{2.3}) is satisfied for any $\mu\in(0,r)$.

For general linear functional differential equations, Lemmas \ref{lem2.1} and \ref{lem2.2} below are known results (see Chapter 3 in \cite{Arino}). However, for convenience
and to emphasize the dependence on the delay kernel $K$, we present Lemmas \ref{lem2.1}
and \ref{lem2.2} here and the proofs are given in Appendix A.
\begin{lemma}\label{lem2.1}
Let $x_{\phi}(t)$ to be the solution of \eqref{2.1} with initial function
$\phi\in BC((-\infty,0],\R).$ Then there exist positive constants $A=A(b,\phi,K)$
and $\gamma=\gamma(a,b)$ such that
\begin{equation}\label{2.6}
|x_{\phi}(t)|\leq A e^{\gamma t}, \quad t\geq 0.
\end{equation}
\end{lemma}

The fundamental solution of the delay differential equation (\ref{2.1}), denoted by $X(t),$
is defined as the solution of \eqref{2.1} with initial condition
\begin{equation*}
X(t)=\left\{ \begin{array}{cc}
1,\quad & t=0,\\ 0,\quad & t<0.
\end{array}\right.
\end{equation*}
Any solutions of (\ref{2.1}) with initial function $\phi\in BC((-\infty,0],\R)$
can be represented through the fundamental solution $X(t)$ as follows.
\begin{lemma}\label{lem2.2}
Let $x_{\phi}(t)$ to be the solution of \eqref{2.1} with initial function
$\phi\in BC((-\infty,0],\R).$ Then
\begin{equation}\label{eq:2.7'}
x_{\phi}(t)= X(t)\phi(0)+b\int_{0}^{t}X(t-s)
\int_{s}^{+\infty}K(\theta)\phi(s-\theta)d\theta ds, \quad t\geq 0.
\end{equation}
\end{lemma}

Properties of the fundamental solution $X(t)$ are closely related to the characteristic function of \eqref{2.1}
defined below.  For any function $f(t):[0,+\infty)\to \R$ which is measurable and satisfies
$$
\left|f(t)\right|\leq a_{1}e^{a_{2}t}, \quad t\in[0,+\infty),
$$
for some constants $a_{1}, a_{2}$, the Laplace transform
\begin{equation*}
\mathcal{L}(f)(\lambda)=\dint_{0}^{+\infty}e^{-\lambda t}f(t)dt, \quad \lambda\in \C
\end{equation*}
exists and is an analytic function of $\lambda$ for $\Re(\lambda)>a_{2}$. Through the Laplace transform $\mathcal{L}(K)$ of the delay kernel $K$, the characteristic function of (\ref{2.1}) is given by
\begin{equation}
h(\lambda)=
\lambda-a-b\mathcal{L}(K)(\lambda).
\end{equation}
It is easy to see that $h(\lambda)$ is well defined and
analytic when $\Re(\lambda) \geq -\mu$, and
\begin{equation}
\mathcal{L}(X)(\lambda) = 1/h(\lambda).
\end{equation}

Now, we can obtain the precise exponential bound of the fundamental solution $X(t)$
in terms of the supremum of the real parts of all roots of the characteristic function $h(\lambda)$.

First, we note that $h(\lambda)$ is analytic when $\Re(\lambda)>-\mu$, and therefore all zeros of  $h(\lambda)$
are isolated. Following the discussion in \cite[Lemma 4.1 in Chapter 1]{Hale93} and \eqref{2.3}, there is a real
number $\alpha_0$ such that all roots of $h(\lambda)=0$ satisfy $\Re(\lambda) < \alpha_0$. Thus,
$\alpha_0 = \sup\{\Re(\lambda): h(\lambda)=0\}$ is well defined. Furthermore, there are only a finite
number of roots in any close subset in the complex plane.

\begin{theorem}
\label{Thm2.3}
Let $\alpha_{0}=\sup\{\Re(\lambda): h(\lambda)=0,\lambda \in\C\}.$
Then
\begin{enumerate}
\item for any $\alpha>\alpha_{0},$ there exists a positive constant $C_{1}=C_{1}(\alpha)$
such that the fundamental solution $X(t)$ satisfies
\begin{equation}\label{2.9}
|X(t)|\leq C_{1}e^{\alpha t},\quad t\geq 0;
\end{equation}
\item for any $\alpha_1 < \alpha_0$, there exist $\bar{\alpha}\in (\alpha_1,\alpha_{0})$ and
a subset $U\subset \mathbb{R}^+$ with measure $m(U)=+\infty$ such that the fundamental solution $X(t)$ satisfies
\begin{equation}\label{2.10}
|X(t)|\geq e^{\bar{\alpha} t},\quad \forall t\in U.
\end{equation}
\end{enumerate}
\end{theorem}

\Proof
1. The proof of (\ref{2.9}) is the same as that of \cite[Theorem 5.2 in Chapter 1]{Hale93} and is omitted here.

2. Let $\alpha_1 < \alpha_0$. Since all zeros of $h(\lambda)$ are isolated, we can take $\bar{\alpha}\in (\alpha_{1},\alpha_{0})$
such that the line $\Re(\lambda)=\bar{\alpha}$ contains no root of the characteristic
equation $h(\lambda)=0$. Next, choose $c>\alpha_{0}$, then
\begin{equation}
\label{eq:Xi}
X(t)=\dfrac{1}{2\pi i}\lim_{T\to+\infty}\int_{c-iT}^{c+iT}\frac{e^{\lambda t}}{h(\lambda)}d\lambda.
\end{equation}

To calculate the integral in \eqref{eq:Xi}, we consider the integration of the function $e^{\lambda t}/h(\lambda)$ around the bounder of the box in the complex plane with boundary $\Gamma = L_{1}M_{1}L_{2}M_{2}$ in the anticlockwise direction, where the segment $L_{1}$ is the set
$\{c+i\tau: -T\leq \tau \leq T\},$ the segment $L_{2}$ is the set $\{\bar{\alpha}+i\tau: -T\leq \tau \leq T\},$
the segment $M_{1}$ is the set $\{u+iT: \bar{\alpha}\leq u \leq c\}$ and
the segment $M_{2}$ is the set $\{u-iT: \bar{\alpha}\leq u \leq c\}$.
From the Cauchy theorem of residues, we obtain
\begin{equation}
\oint_{\Gamma}\frac{e^{\lambda t}}{h(\lambda)}d\lambda = 2\pi i\sum_{j=1}^{m}
\Res_{\lambda=\lambda_{j}}\frac{e^{\lambda t}}{h(\lambda)},
\end{equation}
where $\lambda_{1}, \lambda_{2}, \cdots, \lambda_{m}$ are all roots of $h(\lambda)=0$ inside $\Gamma$ ($m\geq 1$ from the definition of $\alpha_0$, and $m<+\infty$ since $h(\lambda)$ is an analytic function). Here we assume further
$$\bar{\alpha} < \Re(\lambda_1)\leq \Re(\lambda_2)\leq \cdots \leq \Re(\lambda_m)\leq \alpha_0.$$
We note that
\begin{equation*}
\Res_{\lambda=\lambda_{j}}\frac{e^{\lambda t}}{h(\lambda)} = P_{j}(t)e^{\lambda_{j}t},
\end{equation*}
where $P_{j}(t)$ is a nonzero polynomial of $t$ with degree given by the multiplicity
of $\lambda_{j}-1$.
Thus,
\begin{equation}
\label{2.7}\oint_{\Gamma}\frac{e^{\lambda t}}{h(\lambda)}d\lambda = 2\pi i\sum_{j=1}^{m}
P_{j}(t)e^{\lambda_{j}t}.
\end{equation}

There exists a positive constant $\overline{T}$ such that for
$T>\overline{T}$ and $u\geq -\mu$,
\begin{equation*}
\dfrac{\left|h(u+iT)\right|}{T}\geq \dfrac{\sqrt{u^2+T^2}-|a|-|b|\rho}{T} \geq \dfrac{1}{2}.
\end{equation*}
Therefore, when $T$ is large enough,
\begin{eqnarray*}
\left|\int_{M_{1}}\frac{e^{\lambda t}}{h(\lambda)}d\lambda\right|
&=&\left|\int_{\bar{\alpha}+iT}^{c+iT}\frac{e^{\lambda t}}{h(\lambda)}d\lambda\right|\\
&\leq &\int_{\bar{\alpha}}^{c}\dfrac{e^{u t}}{\sqrt{u^2+T^2}-|a|-|b|\rho}du\\
&\leq &\dfrac{2}{T}e^{c t}(c-\bar{\alpha})\to 0 \quad(\mathrm{as}\ T\to +\infty).
\end{eqnarray*}
Similarly, we have
\begin{equation*}
\left|\int_{M_{2}}\frac{e^{\lambda t}}{h(\lambda)}d\lambda\right| \to 0 \quad (\mathrm{as}\ T\to +\infty).
\end{equation*}
Therefore from (\ref{2.7}) we get
\begin{equation*}
\dfrac{1}{2\pi i}\lim_{T\to+\infty}\int_{c-iT}^{c+iT}\frac{e^{\lambda t}}{h(\lambda)}d\lambda
+\dfrac{1}{2\pi i}\lim_{T\to+\infty}\int_{\bar{\alpha}+iT}^{\bar{\alpha}-iT}\frac{e^{\lambda t}}{h(\lambda)}d\lambda
=\sum_{j=1}^{m}P_{j}(t)e^{\lambda_{j}t},
\end{equation*}
which implies
\begin{equation*}
X(t)=X_{\bar{\alpha}}(t)+\sum_{j=1}^{m}P_{j}(t)e^{\lambda_{j}t},
\end{equation*}
where
\begin{equation*}
X_{\bar{\alpha}}(t)=
\dfrac{1}{2\pi i}\lim_{T\to+\infty}\int_{\bar{\alpha}-iT}^{\bar{\alpha}+iT}\frac{e^{\lambda t}}{h(\lambda)}d\lambda.
\end{equation*}

Similar to the proof of (\ref{2.9}), there exists a positive constant
$\bar{C}_1=\bar{C}_1(\bar{\alpha})$ such that $X_{\bar{\alpha}}(t)$ satisfies
\begin{equation}
\left|X_{\bar{\alpha}}(t)\right|\leq \bar{C}_1 e^{\bar{\alpha} t} \;(t\geq 0).\label{2.8}
\end{equation}
Thus
\begin{eqnarray*}
\left|X(t)\right| &\geq&  \left|\sum_{j=1}^{m}P_{j}(t)e^{\lambda_{j}t}\right|
-\left|X_{\bar{\alpha}}(t)\right|
\ \geq\ \left|\sum_{j=1}^{m}P_{j}(t)e^{\lambda_{j}t}\right| - \bar{C}_1e^{\bar{\alpha}t} \\
&=& e^{\bar{\alpha}t} (e^{(\Re(\lambda_1) - \bar{\alpha})t} f(t) - \bar{C}_1),
\end{eqnarray*}
where $f(t) = |\dsum_{j=1}^m P_j(t) e^{(\lambda_j - \Re(\lambda_1))t}|$.

Let $\lambda_j = \beta_j + i \omega_j$, and assume $k$ such that $\beta_j < \beta_m$ when $1\leq j \leq k$, and $\beta_j= \beta_m$ when $k+1\leq j \leq m $, then
\begin{eqnarray*}
f(t) &=& e^{(\beta_m - \beta_1)t}\left|\sum_{j=1}^k e^{-(\beta_m-\beta_j)t} P_j(t) e^{i \omega_j t} + \sum_{j=k+1}^m P_j(t) e^{i \omega_j t} \right|\\
&\geq&\left|\sum_{j=k+1}^m \Re(P_j(t)e^{i\omega_jt})\right| - \sum_{j=1}^k e^{-(\beta_m - \beta_j)t} |P_j(t)|
\end{eqnarray*}
Since $P_j(t)$ are nonzero polynomials, there are a positive constant $\varepsilon > 0$ and a subset $U\subset \mathbb{R}^+$ with measure $m(U)=+\infty$
such that for any $t\in U$,\footnote{It is easy to see that $\sum_{j=k+1}^m\Re(P_j(t)e^{i\omega_j t}) = t^n [\sum_{j=k+1}^m (a_j \cos(\omega_j t) +
b_j\sin(\omega_j t)) + O(t^{-1})]\ (t\to+\infty)$ where $a_j,b_j$ are constants, and $n$ is the highest degree of the polynomials $P_j(t)(j=k+1,\cdots,
m)$. Thus, we can always find a subset $U_0$ with measure $m(U_0)= +\infty$ so that all functions $a_j \cos(\omega_j t) + b_j\sin(\omega_j t) >
\varepsilon\ (k+1\leq j\leq m, \forall t\in U_0)$ for some small positive constant $\varepsilon$ (detail proof is omitted here), and therefore the 
subset $U$ is always possible by taking $U = U_0 \cap (t_0, +\infty)$ with $t_0$ large enough.}

$$\left|\sum_{j=k+1}^m \Re(P_j(t)e^{i\omega_jt})\right|>2\varepsilon.$$
Moreover, since $\dsum_{j=1}^k e^{-(\beta_m - \beta_j)t} |P_j(t)|\to0$ as $t\to +\infty$ and $\Re(\lambda_1) - \bar{\alpha} > 0$, we can further take $U$ such that
$$e^{(\Re(\lambda_1) - \bar{\alpha})t} f(t) - \bar{C}_1 > 1,\quad \forall t\in U$$
and hence (\ref{2.10}) is concluded.
\qed

From  Lemma \ref{lem2.2} and Theorem \ref{Thm2.3}, asymptotical behaviors of solutions $x_\phi(t)$ of equation \eqref{2.1}
are determined by $\alpha_{0}$.

\begin{theorem}\label{Thm2.4}
Let $\alpha_{0}$ be defined as in Theorem \ref{Thm2.3}. Then for any $\alpha>\max\{\alpha_{0},-\mu\}$ there exists
a positive constant $K_{1}=K_{1}(\alpha,\mu)$ such that
\begin{equation}\label{2.16}
|x_{\phi}(t)|\leq K_{1}\|\phi\| e^{\alpha t},\quad t\geq 0,
\end{equation}
where $\mu$ is defined as in \eqref{2.3}. Therefore the trivial solution of \eqref{2.1}
is locally asymptotically stable if and only if $\alpha_{0}<0.$
\end{theorem}

\Proof For any initial function $\phi\in BC((-\infty,0],\R)$
\begin{eqnarray*}
\left|\int_{s}^{+\infty}K(\theta)\phi(s-\theta)d\theta\right|
&\leq &e^{-\mu s}\int_{s}^{+\infty}e^{\mu \theta}K(\theta) \left|\phi(s-\theta)\right|d\theta\\
&\leq &\|\phi\| e^{-\mu s}\int_{s}^{+\infty}e^{\mu \theta}K(\theta)d\theta\\
&\leq & \rho \|\phi\| e^{-\mu s}.
\end{eqnarray*}
Thus from \eqref{eq:2.7'} and Theorem \ref{Thm2.3}, for any $\alpha>\alpha_{0},$
\begin{eqnarray*}
|x_{\phi}(t)|&\leq & |X(t)|\|\phi\|+\int^{t}_{0}|X(t-s)|
\left|\int_{s}^{+\infty}K(\theta)\phi(s-\theta)d\theta\right|ds\\
&\leq & C_{1}\|\phi\|e^{\alpha t}+C_{1}\int^{t}_{0}e^{\alpha (t-s)}
\rho \|\phi\| e^{-\mu s}ds\\
&=&C_{1}\|\phi\|e^{\alpha t}+C_{1}\rho \|\phi\|e^{\alpha t}\int^{t}_{0}e^{-(\alpha+\mu)s}ds\\
&\leq & C_{1}\|\phi\|e^{\alpha t}+C_{1}\rho \|\phi\|\dfrac{e^{-\mu t}+e^{\alpha t}}{|\alpha+\mu|}\\
&=& C_{1}\left(1+\dfrac{2\rho}{|\alpha+\mu|}\right)\|\phi\| e^{\alpha t},
\end{eqnarray*}
Thus, \eqref{2.16} is concluded with
\begin{equation}\label{2.17}
K_{1}(\alpha,\mu)=C_{1}\left(1+\dfrac{2\rho}{|\alpha+\mu|}\right).
\end{equation}
The Theorem is proved.
\qed

For a general distribution density functions $K$, it is not straightforward to obtain sufficient and necessary conditions for $\alpha_0<0$. A sufficient condition is given below.

\begin{theorem}\label{Thm 2.5}
If $(a,b)\in S:=S_{1}\bigcup S_{2},$ where
\begin{eqnarray*}
S_{1}&=&\{(a,b)\in \R^{2}: a<0, a<b<-a \},\\
S_{2}&=&\{(a,b)\in \R^{2}: -(\int_{0}^{+\infty}t K(t)dt)^{-1}<b \leq a\leq 0\},
\end{eqnarray*}
then $\alpha_{0}<0.$
\end{theorem}

\Proof Let $\lambda=\alpha+i\beta\;(\beta>0)$ be a solution of $h(\lambda)=0$. Separating
the real and the imaginary parts, we have
\begin{eqnarray}\label{2.18}
\left\{\begin{array}{ll}
\alpha-a-b\dint_{0}^{+\infty}e^{-\alpha t}K(t)\cos(\beta t)dt=0,\\[0.1cm]
\beta+b\dint_{0}^{+\infty}e^{-\alpha t}K(t)\sin(\beta t)dt=0.
\end{array}\right.
\end{eqnarray}

If $\alpha \geq 0$ and $(a,b)\in S_{1}$, then when $0<b<-a$
\begin{eqnarray*}
\alpha-a-b\int_{0}^{+\infty}e^{-\alpha t}K(t)\cos(\beta t)dt
& > & \alpha+b-b\int_{0}^{+\infty}K(t)dt\\
&= & \alpha+b-b=\alpha \geq 0,
\end{eqnarray*}
and when $a<b\leq 0$
\begin{eqnarray*}
\alpha-a-b\int_{0}^{+\infty}e^{-\alpha t}K(t)\cos(\beta t)dt
&\geq & \alpha-a-|b|\int_{0}^{+\infty}K(t)dt\\
&=& \alpha-(a+|b|)> \alpha \geq 0.
\end{eqnarray*}
Hence, for $(a,b)\in S_{1}$, all
roots of $h(\lambda)=0$ must have negative real parts, i.e., $\alpha_{0}<0$.

If $\alpha \geq 0$ and $(a,b)\in S_{2}$, then for $\beta>0$ we obtain
\begin{eqnarray*}
\beta+b\int_{0}^{+\infty}e^{-\alpha t}K(t)\sin(\beta t)dt
&\geq & \beta+b\int_{0}^{+\infty}K(t)\beta tdt\\
&= &\beta\left(1+b\int_{0}^{+\infty}t K(t)dt\right)>0,
\end{eqnarray*}
which implies that the second equation of (\ref{2.18}) is not satisfied. Hence $\alpha_{0}<0$ when $(a,b)\in S_2$.
\qed

Now we give some properties of the fundamental solution $X(t)$ which are useful for our
estimation of the second moment in the next section.

Obviously, both $X^{2}(t)$ and $X_{s}(t)X_{l}(t)$ have Laplace transforms (here $X_s(t)=X(t-s)$). When $\alpha_0<0$,
the explicit expression of the Laplace transform $\mathcal{L}(X^2)$ is obtained below.

Since $\mathcal{L}(X)=1/h(\lambda)$ and $\alpha_{0}<0$, we have
\begin{equation}\label{2.19}
X(t) = \dfrac{1}{2 \pi}\int_{-\infty}^{\infty} \dfrac{e^{i\omega t}}{h(i\omega)} d\omega.
\end{equation}
Therefore, we obtain
\begin{eqnarray} \label{2.20}
\mathcal{L}(X^2)(\lambda)&=&\int_0^\infty e^{-\lambda t} X^{2}(t)dt
=\dfrac{1}{2\pi}\int_{-\infty}^{\infty} \dfrac{1}{h(i\omega)} \int_{0}^{\infty}
e^{-(\lambda -i\omega)t}X(t) d t d\omega\nonumber\\
&=&\dfrac{1}{2\pi}\int_{-\infty}^{\infty}\dfrac{1}{h(i\omega)h(\lambda-i\omega)} d \omega.
\end{eqnarray}

Let
\begin{equation}\label{2.21}
g(\lambda, s, l) = \frac{\mathcal{L}(X_sX_l)(\lambda)}{\mathcal{L}(X^2)(\lambda)}.
\end{equation}
The function $g(\lambda, s,l)$ is crucial for the characteristic function of \eqref{1.1}.
Similar to the above argument, we obtain an explicit expression of $g(\lambda, s,l)$ given in Lemma \ref{lem2.6}.
\begin{lemma}
\label{lem2.6}
Let $g(\lambda, s,l)$ defined as in \eqref{2.21}, then
\begin{equation}\label{2.22}
g(\lambda, s,l) =\left\{
\begin{array}{ll}
\dfrac{\int_{-\infty}^{+\infty}\dfrac{e^{-\lambda s}e^{i\omega(s-l)}}{h(i\omega)h(\lambda-i\omega)}d\omega }
{\int_{-\infty}^{+\infty}\dfrac{1}{h(i\omega)h(\lambda-i\omega)}d\omega}, &\, s\geq l>0,\\
\\
\dfrac{\int_{-\infty}^{+\infty}
\dfrac{e^{-\lambda l}e^{i\omega(l-s)}}{h(i\omega)h(\lambda-i\omega)}d\omega }
{\int_{-\infty}^{+\infty}
\dfrac{1}{h(i\omega)h(\lambda-i\omega)}d\omega}, &\, 0\leq s<l.
\end{array}\right.
\end{equation}
\end{lemma}

\Proof From (\ref{2.19}),  we have
\begin{eqnarray}
\mathcal{L}(X_{s}X_{l})&=&\int_{0}^{+\infty}e^{-\lambda t}X(t-s)X(t-l)dt\nonumber\\
&=& \int_{0}^{+\infty}e^{-\lambda t}\dfrac{1}{2 \pi}\int_{-\infty}^{\infty}\dfrac{e^{i\omega(t-l)}}
{h(i\omega)} d\omega X(t-s)dt \nonumber\\
&=& \dfrac{1}{2 \pi}\int_{-\infty}^{\infty}\dfrac{e^{-i\omega l}e^{-(\lambda-i\omega) s}}
{h(i\omega)}\int_{0}^{+\infty}e^{-(\lambda-i\omega)(t-s)}X(t-s)dt d\omega \nonumber\\
&=&\left\{\begin{array}{ll} \dfrac{1}{2 \pi}\dint_{-\infty}^{\infty}\dfrac{e^{-\lambda s}e^{i\omega (s-l)}}
{h(i\omega)h(\lambda-i\omega)}d\omega\;(s\geq l\geq 0) \\
\\
\dfrac{1}{2 \pi}\dint_{-\infty}^{\infty}\dfrac{e^{-\lambda l}e^{i\omega (l-s)}}
{h(i\omega)h(\lambda-i\omega)}d\omega \;(0\leq s\leq l ).
\end{array}\right.\label{2.23}
\end{eqnarray}
Thus \eqref{2.22} is followed from (\ref{2.20}), (\ref{2.21}) and (\ref{2.23}).
\qed

From \eqref{2.22}, we have
\begin{equation}\label{2.24}
g(\lambda, s,0) =\dfrac{\int_{-\infty}^{+\infty}
\dfrac{e^{-i\omega s}}{h(i\omega)h(\lambda-i\omega)}d\omega }
{\int_{-\infty}^{+\infty}
\dfrac{1}{h(i\omega)h(\lambda-i\omega)}d\omega}.
\end{equation}

The following Lemma gives an important estimation of $g(\lambda, s,l)$, with the proof given in Appendix B.
\begin{lemma}
\label{lem2.7}
Let $g(\lambda, s,l)$ defined as in \eqref{2.21}. Then
when  $\Re(\lambda)>\max\{2\alpha_0, a, a+\alpha_0, -\mu\}$, for any $\varepsilon >0,$ there exists a constant
$T_{0}=T_{0}(\varepsilon)$ independent to $s$ and $l$ such that for $|\lambda|>T_{0}$,
\begin{equation}\label{2.25}
\left|g(\lambda, s,l)\right|\leq
\left\{\left.\begin{array}{cc}
\dfrac{e^{-(T_{0}-a)l}e^{-as}}{1-\varepsilon}
+\dfrac{\varepsilon e^{-T_{0}s}}{1-\varepsilon},& (s\geq l\geq 0), \\[0.2cm]
\dfrac{e^{-(T_{0}-a)s}e^{-al}}{1-\varepsilon}+\dfrac{\varepsilon e^{-T_{0}l}}{1-\varepsilon},& (0\leq s<l).
\end{array}\right.\right.
\end{equation}
Moreover, for $\Re(\lambda)>\max\{2\alpha_0, a, a+\alpha_0,-\mu\},$
\begin{equation}\label{2.26}
\lim_{|\lambda|\to +\infty}\left|g(\lambda,s,l)\right|=0.
\end{equation}
\end{lemma}

\section{Moment boundedness of the equation with noise perturbation}

Now we consider the equation \eqref{1.1}, i.e., $\sigma_{i} ~(i=0,1,2)$ are not
all zeros. In this section, two main results are obtained: Theorem \ref{Thm3.3} for the sufficient
condition of the exponential stability of the first moment, and Theorem \ref{Thm3.8} for the
characteristic equation that implies the
boundedness criteria for the second moments of solutions of equation \eqref{1.1}.

The existence and uniqueness theorem  for the stochastic differential
delay equations have been established in \cite{Ito,Mao1,Mohammed84}. Hence using the
fundamental solution $X(t)$ in the previous section, the solution $x(t;\phi)$ of (\ref{1.1}) with the
initial function  $\phi\in BC((-\infty,0],\R)$ is a
1-dimensional stochastic process given by It\^{o} integral as follows:
\begin{equation}
\label{3.1}
\begin{array}{rcl}
x(t;\phi)&=&x_{\phi}(t)+\dint^{t}_{0}X(t-s)\Big(\sigma_{0}+\sigma_{1}x(s;\phi)\\
&&{} \qquad\qquad+\sigma_{2}\dint_{0}^{+\infty}K(\theta)x(s-\theta;\phi)d\theta\Big)dW_s,
\end{array}\quad t\geq 0,
\end{equation}
where $x_{\phi}(t)$ is the solution of \eqref{2.1} defined by
\eqref{eq:2.7'} and $W_s$ is a 1-dimensional Wiener process.

The first and second moments of $x(t;\phi)$ are very important for investigating
the behavior of the solutions and are studied in this paper.  Now we state definitions of
the $p^{\mathrm{th}}$ moment exponential stability and the $p^{\mathrm{th}}$ moment boundedness. Here we denote by $E$
the mathematical expectation.

\begin{definition}\label{def3.1}
The solution of  \eqref{1.1} is said to be the first moment exponentially
stable if there exist two positive constants $\gamma$ and $R$ such that
\begin{equation*}
|E(x(t;\phi))|\leq R\|\phi\|e^{-\gamma t},\quad t\geq 0,
\end{equation*}
for all $\phi\in BC((-\infty,0],\R)$. When $p\geq 2,$ the solution of
\eqref{1.1} is said to be the $p^{\mathrm{th}}$ moment exponentially stable if there
exist two positive constants $\gamma$ and $R$ such that
\begin{equation*}
E\big(|x(t;\phi)-E(x(t;\phi))|^{p}\big)\leq R\|\phi\|^{p}e^{-\gamma t},\quad t\geq 0,
\end{equation*}
for all $\phi\in BC((-\infty,0],\R)$.
\end{definition}

\begin{definition}\label{def3.3}
For $p\geq 2,$ the solution of \eqref{1.1} is said to be the $p^{\mathrm{th}}$ moment
bounded if there exists a positive constant $\tilde{R}=\tilde{R}(\|\phi\|^p)$
such that
\begin{equation*}
E\big(|x(t;\phi)-E(x(t;\phi))|^{p}\big)\leq \tilde{R},\quad t\geq 0,
\end{equation*}
for all $\phi\in BC((-\infty,0],\R)$. Otherwise, the $p^{\mathrm{th}}$ moment
is said to be unbounded.
\end{definition}

We first investigate the exponential stability of the first moment.

\subsection{First moment stability}
From \eqref{3.1}, it is easy to have $Ex(t;\phi)=x_{\phi}(t)$ from the It\^{o} integral, and therefore Theorem \ref{Thm2.4} yields the following result.

\begin{theorem}\label{Thm3.3}
Let $\alpha_{0}$ be defined as in Theorem \ref{Thm2.3}.
Then for any $\alpha>\max\{\alpha_{0},-\mu\}$ there exist a constant
$K_{1}=K_{1}(\alpha, \mu)$ defined in \eqref{2.17} such that
\begin{equation}\label{3.3}
|Ex(t;\phi)|\leq K_{1}\|\phi\| e^{\alpha t},
\quad t\geq 0.
\end{equation}
Therefore if $\alpha_{0}<0$ the solution of \eqref{1.1} is first moment exponentially stable.
\end{theorem}

Theorem \ref{Thm3.3} indicates that the stability condition of the first moment is the same as that of
the unperturbed equation \eqref{2.1}. The stability is determined by coefficients
$a$ and $b$ and is independent of the parameters $\sigma_i\ (i=0,1,2)$.

\subsection{Second moment boundedness}
Now we study the behavior of the second moment. Let $x(t;\phi)$ be a solution of (\ref{1.1}), and define
\begin{eqnarray*}
\tilde{x}(t;\phi)&=&x(t;\phi)-Ex(t;\phi),\quad M(t)=E(\tilde{x}^{2}(t;\phi))
\end{eqnarray*}
and
$$
N(t;s,l)=E\big(\tilde{x}(t-s;\phi)\tilde{x}(t-l;\phi)\big)\quad (t, s, l \geq 0).
$$
Then $M(t)=N(t;0,0)$ is the second moment of $x(t;\phi)$. Obviously, when $ t \leq 0$,
$\tilde{x}(t;\phi)=E\tilde{x}(t;\phi)=M(t)=0$, and when $s\geq t$ or $l \geq t$, $N(t;s,l)=0$.

We introduce following notations:
\begin{eqnarray*}
P(t)&=&\left\{
\left.\begin{array}{ll}
\left(\sigma_{0}+\sigma_{1}Ex(t;\phi)
+\sigma_{2}\dint_{0}^{+\infty}K(\theta)Ex(t-\theta;\phi)d\theta \right)^{2}, & t \geq 0,\\[0.4cm]
0, & t<0,
\end{array}\right.\right.\\
Q(t)&=&\left\{
\left.\begin{array}{ll}
\sigma^{2}_{1}M(t)+2\sigma_{1}\sigma_{2}\dint_{0}^{t}K(s)N(t;s,0)ds & \\[0.2cm]
\qquad\qquad {}
+\sigma^{2}_{2}\dint_{0}^{t}\dint_{0}^{t}K(s)K(l)N(t;s,l)dsdl,\quad & t \geq 0,\\[0.3cm]
0, & t<0,
\end{array}\right.\right.\\
F(t)&=&\dint^{t}_{0}X^{2}(t-s)P(s)ds \;(t \geq 0).
\end{eqnarray*}
Applying It\^{o} integral, a tedious calculation
yields
\begin{equation}\label{3.4}
N(t;s,l)=\int_{0}^{(t-s)\wedge(t-l)}X(t-s-\theta)X(t-l-\theta)\left(P(\theta)+Q(\theta)\right)d\theta,
\end{equation}
where $(t-s)\wedge(t-l)=\min\left\{t-s, t-l\right\}$. Therefore
\begin{equation}\label{3.5}
N(t;s,0)=\int_{0}^{t-s}X(t-\theta)X(t-s-\theta)\left(P(\theta)+Q(\theta)\right)d\theta
\end{equation}
and
\begin{equation}\label{3.6}
M(t)=\int^{t}_{0}X^{2}(t-\theta)\left(P(\theta)+Q(\theta)\right)d\theta.
\end{equation}

\subsubsection{Additive noise}
When $\sigma_{1}=\sigma_{2}=0$, we have only additive noise
and the second moment becomes
\begin{equation}
\label{2.15}
M(t)=\sigma^{2}_{0}\int^{t}_{0}X^{2}(s)ds.
\end{equation}
In this case, from Theorem 2.3, the sufficient conditions for the second moment $M(t)$
to be bounded or unbounded are given as follows.

\begin{theorem}\label{Thm3.4}
Let $\alpha_{0}$ be defined as in Theorem \ref{Thm2.3}.
When $\sigma_{1}=\sigma_{2}=0,$ then
\begin{enumerate}
\item if $\alpha_{0}<0$, the second moment of
\eqref{3.1} is bounded. Moreover, for any $\alpha\in (\alpha_{0},0),$ there exists
a constant $C_{1}=C_{1}(\alpha)$ (as in Theorem 2.3) such that
\begin{equation*}
\left| M(t)-M_{1}\right| \leq \dfrac{C_{1}^2\sigma^{2}_{0}e^{2\alpha t}}
{|2\alpha|},\quad t\geq 0,
\end{equation*}
where
\begin{equation*}
M_{1}=\lim_{t\to +\infty}M(t)=\sigma^{2}_{0}\int^{+\infty}_{0}X^{2}(s)ds
\leq \dfrac{C_{1}^2\sigma^{2}_{0}}{|2\alpha|};
\end{equation*}
\item if $\alpha_{0}>0$, the second moment of \eqref{3.1} is unbounded.
\end{enumerate}
\end{theorem}

\Proof1. If $\alpha_{0}<0$, for any $\alpha \in(\alpha_{0}, 0)$
we obtain
\begin{equation*}
M(t)\leq \sigma^{2}_{0}\int^{t}_{0}C_{1}^2e^{2\alpha s}ds
=\dfrac{C_{1}^2\sigma^{2}_{0}}{|2\alpha|}\left(1-e^{2\alpha t}\right)
\leq \dfrac{C_{1}^2\sigma^{2}_{0}}{|2\alpha|}
\end{equation*}
from Theorem \ref{2.3}, and hence $M(t)$ is bounded.
Moreover,
\begin{equation*}
\left|M(t)-M_{1}\right|=\left|\sigma^{2}_{0}\int^{+\infty}_{t}X^{2}(s)ds\right|
\leq C_{1}^2\sigma^{2}_{0}\int^{+\infty}_{t}e^{2\alpha s}ds
=\dfrac{C_{1}^2\sigma^{2}_{0}e^{2\alpha t}}{|2\alpha|}
\end{equation*}
and
\begin{equation*}
M_{1}\leq C_{1}^2\sigma^{2}_{0}\int^{+\infty}_{0}e^{2\alpha s}ds
=\dfrac{C_{1}^2\sigma^{2}_{0}}{|2\alpha|}.
\end{equation*}

2. If $\alpha_{0} > 0$, from Theorem \ref{Thm2.3}, there exist
$\bar{\alpha}\in (0,\alpha_{0})$ and a closed subset $U\subset \mathbb{R}^+$ with $m(U) = +\infty$ such
that
\begin{equation*}
|X(t)|\geq e^{\bar{\alpha} t},\quad \forall t\in U.
\end{equation*}
Thus from \eqref{2.15},
\begin{equation*}
\lim_{t\to +\infty}M(t)\geq \sigma^{2}_{0}\int_{U}X^{2}(s)ds
\geq \sigma^{2}_{0} \int_{U}e^{2\bar{\alpha}s}ds
\geq \sigma^{2}_{0}  m(U)=+\infty,
\end{equation*}
which implies
\begin{equation*}
\lim_{t\to +\infty}M(t)=+\infty.
\end{equation*}
Therefore, the second moment is unbounded.
\qed

\begin{remark}
The critical case $\alpha_{0}=0$ is not discussed here and the stability issue remains open.
\end{remark}

\subsubsection{General cases ($\sigma_{1},\,\sigma_{2}$ are not all zeros)}
First, we note a very special situation that $\sigma_i\ (i=0,1,2)$ satisfy the following condition:
\vspace{0.2cm}
\begin{center}
\begin{minipage}{10cm}
\begin{description}
\item[\textbf{H}:] $\sigma_0 = 0$, and there is a constant $\lambda$ such that $h(\lambda) = 0$ and $\sigma_1 + \sigma_2 \mathcal{L}(K)(\lambda) = 0$.
\end{description}
\end{minipage}
\end{center}
\vspace{0.2cm}
In this situation, it is easy to verify that $x(t) = e^{\lambda t}\ (t\in \mathbb{R})$ is a solution of \eqref{1.1} with initial function $\phi(\theta) = e^{\lambda \theta}\ (\theta\leq 0)$, and therefore the corresponding second moment $M(t) = 0$. This is a very rare situation to have a deterministic solution of a stochastic delay differential equation, and is excluded in following discussions.

The following result gives the sufficient conditions for the second moment of (\ref{3.1}) to be
unbounded when the trivial solution of \eqref{2.1} is unstable.

\begin{theorem}
\label{Thm3.6}
Let $\alpha_{0}$ be defined as in Theorem \ref{Thm2.3}. If $\alpha_{0}>0$ and the condition \textnormal{\textbf{H}} is not satisfied,
then the second moment of \eqref{3.1} is unbounded.
\end{theorem}
\Proof We only need to show that there is a special solution $x(t;\phi)$ such that the corresponding second moment is unbounded. First, we note
$$Q(t) = E(\sigma_1 \tilde{x}(t) + \sigma_2 \int_0^t K(s)\tilde{x}(t-s) d s)^2\geq 0,$$
and therefore
$$M(t) \geq F(t)=\int_0^t X^2(t-s) P(s) d s.$$

Now, let $\lambda = \alpha + i\beta$ be a solution of $h(\lambda) = 0$ with $0<\alpha\leq \alpha_0$, then $x_{\phi}(t)= \Re(e^{\lambda t})$ is a solution of \eqref{2.1} with initial function $\phi(\theta)=\Re(e^{\lambda \theta})\ (\theta \leq 0)$. Hence, for the solution $x(t;\phi)$ of \eqref{1.1} with this particular initial function, we have
\begin{eqnarray*}
P(t) &=& \left(\Re\left[\sigma_0 + e^{\lambda t}(\sigma_1 + \sigma_2 \mathcal{L}(K)(\lambda))\right]\right)^2\\
&=&\left(\sigma_0 + e^{\alpha t}\Re[e^{i\beta t}(\sigma_1 + \sigma_2 \mathcal{L}(K)(\lambda))]\right)^2.
\end{eqnarray*}
Since the condition \textbf{H} is not satisfied, we have either $\sigma_0\not=0$ or $\sigma_1 + \sigma_2 \mathcal{L}(K)(\lambda)\not=0$. Thus, from $\alpha > 0$, and following the proof of Theorem \ref{Thm2.3}, there are a subset $U\subset \mathbb{R}^+$ with measure $m(U) = +\infty$ and $\varepsilon>0$ such that
$$X^2(t) > e^{\bar{\alpha}t},\quad P(t) > \varepsilon,\quad \forall t\in U,$$
where $0< \bar{\alpha} < \alpha_0$. Thus,
$$\lim_{t\to\infty}F(t) = \lim_{t\to+\infty}\int_0^t X^2(t - s) P(s) d s \geq \varepsilon m(U) = +\infty,$$
which implies that the second moment is unbounded.
\qed

\begin{remark}
\label{remark3.7}
From the proof, for any $\lambda$ with $\Re(\lambda)>0$ such that $h(\lambda)=0$, if either $\sigma_0 \not=0$ or $\sigma_1 + \sigma_2 \mathcal{L}(K)(\lambda)\not=0$, the second moment of the solution of \eqref{1.1} with initial function $\phi(\theta) = \Re(e^{\lambda \theta})\ (\theta \leq 0)$ is unbounded.
\end{remark}

In the following discussions, we always assume $\alpha_{0}< 0$.

Now we study the second moment through the method of Laplace transform. First we note that both $M(t)$ and $N(t;s,l)$
have Laplace transforms, for detail proofs refer to Lemmas \ref{lem3.9} and \ref{lem3.10} below.

The following theorem presents the characteristic function of \eqref{1.1} and establishes the boundedness criteria for  the second moment of the solutions of \eqref{1.1}.

\begin{theorem}\label{Thm3.8}
Let $\alpha_{0}$ be defined as in Theorem \ref{Thm2.3} and assume $\alpha_0 < 0$. Define
\begin{equation}\label{3.7}
H(\lambda)=\lambda-\left(2a+\sigma^{2}_{1}\right)-2\left(b+\sigma_{1}\sigma_{2}\right)f_{1}(\lambda)
-\sigma^{2}_{2}f_{2}(\lambda),
\end{equation}
where
\begin{equation}\label{3.8}
\left.\left.\begin{array}{ll}
f_{1}(\lambda)=\dint_{0}^{+\infty}K(s)g(\lambda, s, 0)ds, & \\
f_{2}(\lambda)=\dint_{0}^{+\infty}\int_{0}^{+\infty}K(s)K(l)g(\lambda, s,l)ds dl,  &
\end{array}\right.\right.
\end{equation}
and $g(\lambda, s,l)$ is defined as in \eqref{2.21}. Then
\begin{enumerate}
\item
if all roots of the characteristic equation $H(\lambda)=0$ have negative real parts,
the second moment for any solution of \eqref{1.1} is bounded, and approaches a
constant exponentially as $t\rightarrow +\infty$;
\item if the characteristic equation $H(\lambda)=0$ has a root with the positive real part, and the condition \textnormal{\textbf{H}} is not satisfied, the second moment of \eqref{3.1} is unbounded.
\end{enumerate}
\end{theorem}

From Theorem \ref{Thm3.8}, $H(\lambda)$ is the characteristic function for the second moment
boundedness of the stochastic delay differential equation \eqref{1.1}. We note that the
characteristic function is independent of the coefficient $\sigma_0$. But as we can see in the proof below,
when the second moment is bounded, the limit $\dlim_{t\to\infty}M(t)$ depends on $\sigma_0$. To prove Theorem \ref{Thm3.8},
we first give some lemmas.
\begin{lemma}\label{lem3.9}
For any $\alpha\in(\alpha_{0}, 0),$ there exist a positive constant
$K_{2}=K_{2}(\alpha,\phi)$  such that
\begin{equation}\label{3.9}
F(t)\leq K_{2}(1-e^{2\alpha t}),\quad t \geq 0.
\end{equation}
\end{lemma}

\Proof Since $\mu>0,\,\alpha_{0}<0$ and (\ref{2.3}), from Theorem \ref{Thm3.3},
for any $\alpha \in (\alpha_{0},0)$, there exists a positive constant $K_1$ such that
\begin{eqnarray}\label{3.10}
&&\left|\int_{0}^{+\infty}K(\theta)Ex(s-\theta;\phi)d\theta\right| \nonumber\\
&\leq &\left|\int_{0}^{s}K(\theta)Ex(s-\theta;\phi)d\theta\right|
+\left|\int_{s}^{+\infty}K(\theta)Ex(s-\theta;\phi)d\theta\right| \nonumber\\
&\leq & \int_{0}^{s}K(\theta)K_{1}\|\phi\|e^{\alpha(s-\theta)}d\theta
+\int_{s}^{+\infty}K(\theta)|\phi(s-\theta)|d\theta \nonumber \\
&\leq &K_{1}\|\phi\|\int_{0}^{s}K(\theta)d\theta
+\|\phi\|\int_{s}^{+\infty}K(\theta)d\theta\nonumber\\
&\leq &\left(1+K_{1}\right)\|\phi\|.
\end{eqnarray}
Thus from Theorem \ref{Thm2.3}, for any $\alpha\in(\alpha_{0}, 0),$
\begin{eqnarray*}
F(t)&\leq& \int^{t}_{0}C_{1}^{2}e^{2\alpha(t-s)}\left(\sigma_{0}+\sigma_{1}Ex(s;\phi)
+\sigma_{2}\int_{0}^{+\infty}K(s-\theta;\phi)Ex(\theta;\phi)d\theta \right)^{2}ds\\
&\leq&C_{1}^{2}e^{2\alpha t} \int^{t}_{0}e^{-2\alpha s}
\Big(|\sigma_{0}|+|\sigma_{1}|K_{1}\|\phi\|e^{\alpha s}
+|\sigma_{2}|\left(1+K_{1}\right)\|\phi\|\Big)^{2}ds\\
&\leq& C_{1}^{2}e^{2\alpha t} \int^{t}_{0}\Big(|\sigma_{0}|+|\sigma_{1}|K_{1}\|\phi\|
+|\sigma_{2}|\left(1+K_{1}\right)\|\phi\|\Big)^{2}
e^{-2\alpha s}ds\\
&=&C_{1}^{2}\Big(|\sigma_{0}|+|\sigma_{1}|K_{1}\|\phi\|+|\sigma_{2}|
\left(1+K_{1}\right)\|\phi\|\Big)^{2}
\dfrac{1-e^{2\alpha t}}{|2\alpha|}\\
&=&K_{2}(1-e^{2\alpha t}),
\end{eqnarray*}
where
$$
K_{2}(\alpha,\mu)=\dfrac{C_{1}^{2}}{|2\alpha|}\Big(|\sigma_{0}|+|\sigma_{1}|K_{1}\|\phi\|
+|\sigma_{2}|\left(1+K_{1}\right)\|\phi\|\Big)^{2}.
$$
The Lemma is proved.
\qed

A direct consequence of Lemma \ref{lem3.9} is that $F(t)$ has Laplace transform.
In the following, we have similar estimation for $M(t)$.

\begin{lemma}\label{lem3.10}
Let $\alpha_0$ be defined as in Theorem \ref{Thm2.3} and assume $\alpha_0 < 0$, then
\begin{equation}\label{3.11}
M(t)\leq K_2 e^{C_1^2(|\sigma_1| + |\sigma_2|^2)t},\quad t\geq 0.
\end{equation}
\end{lemma}

\Proof From \eqref{3.6}, we have
\begin{equation}\label{3.12}
M(t) = F(t) + \int_0^t X^2(t-s)Q(s) d s.
\end{equation}
To estimate the integral, we note that
\begin{eqnarray}\label{3.13}
\left|N(t;s,l)\right|&=&\left|E\left(\tilde{x}(t-s)\tilde{x}(t-l)\right)\right|
\leq  \left(E(\tilde{x}^2(t-s))\right)^{\frac{1}{2}}
\left(E(\tilde{x}^2(t-l))\right)^{\frac{1}{2}} \nonumber\\
&\leq &\dfrac{M(t-s)+M(t-l)}{2}
\end{eqnarray}
by the Cauchy-Schwarz inequality.
Therefore, we have
\begin{eqnarray}\label{3.14}
\left|\int_{0}^{t}K(s)N(t;s,0)ds\right|
&\leq &\int_{0}^{t}K(s)\left|N(t;s,0)\right|ds \nonumber\\
&\leq &\dfrac{1}{2}M(t)\int_{0}^{t}K(s)ds+\dfrac{1}{2}\int_{0}^{t}K(s)M(t-s)ds \nonumber\\
&\leq& \dfrac{1}{2}M(t)+\dfrac{1}{2}\int_{0}^{t}K(s)M(t-s)ds
\end{eqnarray}
and
\begin{eqnarray}\label{3.15}
&&\left|\int_{0}^{t}\int_{0}^{t}K(s)K(l)N(t;s,l)ds dl \right|\nonumber\\
&\leq &\int_{0}^{t}\int_{0}^{t}K(s)K(l)\left|N(t;s,l)\right|ds dl\nonumber\\
&\leq & \int_{0}^{t}\int_{0}^{t}K(s)K(l)\dfrac{M(t-s)+M(t-l)}{2}ds dl \nonumber\\
&=&\dfrac{1}{2}\int_{0}^{t}K(s)M(t-s)ds\int_{0}^{t}K(l)dl
+\dfrac{1}{2}\int_{0}^{t}K(l)M(t-l)dl\int_{0}^{t}K(s)ds \nonumber\\
&\leq & \dfrac{1}{2}\int_{0}^{t}K(s)M(t-s)ds
+\dfrac{1}{2}\int_{0}^{t}K(l)M(t-l)dl \nonumber\\
&=&\int_{0}^{t}K(s)M(t-s)ds.
\end{eqnarray}
It is easy to verify that $M(t)$ is increasing on $[0,+\infty)$, and from \eqref{2.2},
\begin{equation}\label{3.16}
\int_{0}^{t}K(s)M(t-s)ds\leq M(t).
\end{equation}
Thus, from Lemma \ref{lem3.9} and (\ref{3.12})-(\ref{3.16}), for any
$\alpha\in(\alpha_{0}, 0),$
\begin{eqnarray*}
M(t)&\leq & K_{2}+\int_{0}^{t}X^2(t-s)\left(\sigma^{2}_{1}M(s)+2|\sigma_{1}\sigma_{2}|M(s)
+\sigma^{2}_{2}M(s)\right)ds\\
&\leq & K_{2}+C_{1}^2\left(|\sigma_{1}|+|\sigma_{2}|\right)^{2}\int_{0}^{t}M(s)ds.
\end{eqnarray*}
Finally, applying the Gronwall inequality, we obtain
\begin{equation*}
M(t)\leq  K_{2}e^{C_{1}^2\left(|\sigma_{1}|+|\sigma_{2}|\right)^{2}t},
\end{equation*}
and \eqref{3.11} is proved.
\qed

Lemma \ref{lem3.10} indicates that  $M(t)$ has Laplace transform. Furthermore, from
(\ref{3.11}) and  (\ref{3.13}), $N(t;s,l)~(0\leq s, l\leq t)$  also has
Laplace transform.

\begin{lemma}\label{lem3.11}
Let $Q(t)$ and $M(t)$ defined as previous. Then
\begin{equation}\label{3.17}
\mathcal{L}(Q)(\lambda) = (\sigma_1^2  + 2 \sigma_1 \sigma_2 f_1(\lambda)
+ \sigma_2^2 f_2(\lambda)) \mathcal{L}(M)(\lambda).
\end{equation}
\end{lemma}
\Proof First, from the expression of $Q(t)$, we have for $t\geq 0,$
\begin{eqnarray}\label{3.18}
\mathcal{L}(Q)(\lambda)&=&\sigma^{2}_{1}\mathcal{L}(M)(\lambda)
+2\sigma_{1}\sigma_{2}\int_{0}^{+\infty}e^{-\lambda t}\int_{0}^{t}K(s)N(t;s,0)ds dt  \nonumber\\
&&\quad{}+\sigma^{2}_{2}\int_{0}^{+\infty}e^{-\lambda t}\int_{0}^{t}\int_{0}^{t}K(s)K(l)N(t;s,l)ds dl dt.
\end{eqnarray}
A direct calculation yields
\begin{eqnarray}\label{3.19}
&&\int_{0}^{+\infty}e^{-\lambda t}\int_{0}^{t}\int_{0}^{t}K(s)K(l)N(t;s,l)ds dl dt \nonumber\\
&=& \int_{0}^{+\infty}\int_{l}^{+\infty}e^{-\lambda t}\int_{0}^{t}K(s)K(l)N(t;s,l)ds dt dl  \nonumber \\
&=& \int_{0}^{+\infty}K(l)\int_{0}^{l}K(s)\int_{l}^{+\infty}e^{-\lambda t}N(t;s,l) dt ds dl \nonumber\\
&&{}
+\int_{0}^{+\infty}K(l)\int_{l}^{+\infty}K(s)\int_{s}^{+\infty}e^{-\lambda t}N(t;s,l)dtds  dl \nonumber \\
&=& \int_{0}^{+\infty}K(l)\int_{0}^{l}K(s)\int_{0}^{+\infty}e^{-\lambda t}N(t;s,l)dt ds dl \nonumber\\
&&{}
+\int_{0}^{+\infty}K(l)\int_{l}^{+\infty}K(s)\int_{0}^{+\infty}e^{-\lambda t}N(t;s,l)dtds  dl \nonumber \\
&=& \int_{0}^{+\infty}K(l)\int_{0}^{+\infty}K(s)\mathcal{L}\left(N(t;s,l)\right)ds  dl,
\end{eqnarray}
and similarly
\begin{equation}\label{3.20}
\int_{0}^{+\infty}e^{-\lambda t}\int_{0}^{t}K(s)N(t;s,0)ds dt
= \int_{0}^{+\infty}K(s)\mathcal{L}(N(t;s,0)) ds.
\end{equation}
Since
\begin{equation*}
N(t;s,l)=\int_{0}^{(t-s)\wedge(t-l)}X(t-s-\theta)X(t-l-\theta)\left(P(\theta)+Q(\theta)\right)d\theta,
\end{equation*}
we have
\begin{equation}\label{3.21}
\mathcal{L}(N(t;s,l))
=\mathcal{L}(X_{s}X_{l})\left(\mathcal{L}(P)+\mathcal{L}(Q)\right).
\end{equation}
We note
\begin{equation}\label{3.22}
\mathcal{L}(M)=\mathcal{L}(X^2)\left(\mathcal{L}(P)+\mathcal{L}(Q)\right)
\end{equation}
by applying the Laplace transform to both sides of \eqref{3.6}. Therefore, for any $s, l \in[0, t]$,
\eqref{3.21} and \eqref{3.22} yield
\begin{equation*}
\mathcal{L}(N(t;s,l)) = \dfrac{\mathcal{L}(X_{s}X_{l})}{\mathcal{L}(X^{2})}\mathcal{L}(M)
=g(\lambda,s,l)\mathcal{L}(M)
\end{equation*}
and
\begin{equation*}
\mathcal{L}(N(t;s,0)) = \dfrac{\mathcal{L}(XX_{s})}{\mathcal{L}(X^{2})}\mathcal{L}(M)
=g(\lambda,s,0)\mathcal{L}(M).
\end{equation*}

Thus, from \eqref{3.19} and \eqref{3.20}, we obtain
\begin{equation}\label{3.23}
\int_{0}^{+\infty}e^{-\lambda t}\int_{0}^{t}\int_{0}^{t}K(s)K(l)N(t;s,l)ds dl dt =f_{2}(\lambda)\mathcal{L}(M)
\end{equation}
and
\begin{equation} \label{3.24}
\int_{0}^{+\infty}e^{-\lambda t}\int_{0}^{t}K(s)N(t;s,0)ds dt
= f_{1}(\lambda)\mathcal{L}(M).
\end{equation}
Finally, \eqref{3.17} is concluded from (\ref{3.18}), (\ref{3.23}) and  (\ref{3.24}).
\qed

Now, we are ready to prove Theorem \ref{Thm3.8}.

\noindent\textit{Proof of Theorem \ref{Thm3.8}.} From \eqref{3.17} and \eqref{3.22}, we obtain
\begin{equation}\label{3.25}
\mathcal{L}(M)=\frac{1}{\mathcal{L}(X^2)^{-1} - (\sigma_1^2 + 2 \sigma_1 \sigma_2 f_1(\lambda)
+ \sigma_2^2 f_2(\lambda))}\mathcal{L}(P).
\end{equation}

To obtain $\mathcal{L}(X^2)^{-1}$, multiplying  $2X(t)$ to both sides of (\ref{2.1}), we have
\begin{equation}\label{3.26}
\dfrac{dX^{2}(t)}{dt}=2aX^{2}(t)+2bX(t)\int_{0}^{+\infty}K(s)X(t-s)ds.
\end{equation}
Taking the Laplace transform to both sides of (\ref{3.26}) yields
\begin{equation*}
-1+\lambda\mathcal{L}(X^2)=2a\mathcal{L}(X^2)+2b\int_{0}^{+\infty}K(s)\mathcal{L}(XX_{s})ds,
\end{equation*}
which gives
\begin{equation}\label{3.27}
\dfrac{1}{\mathcal{L}(X^2)}=\lambda-2a-2bf_{1}(\lambda).
\end{equation}
Now, from (\ref{3.25}) and (\ref{3.27}) we obtain
\begin{equation}\label{3.28}
\mathcal{L}(M)=\dfrac{1}{H(\lambda)}\mathcal{L}(P).
\end{equation}
Let $Y(t)=\mathcal{L}^{-1}\left(H^{-1}\left(\lambda\right)\right)$, then \eqref{3.28} yields
\begin{equation}\label{3.29}
M(t)=Y(t)\ast P(t)=\int_{0}^{t}Y(s)P(t-s)ds,
\end{equation}
where $\ast$ denotes the  convolution product. Now, the exponential bounds for $M(t)$ can be obtained from \eqref{3.29}.

From Lemma \ref{lem2.7} and note $\alpha_0<0$,
\begin{equation*}
\lim_{\Re(\lambda)\to +\infty} |f_1(\lambda)|  = \lim_{\Re(\lambda)\to +\infty} |f_2(\lambda)| = 0.
\end{equation*}
Furthermore, $H(\lambda)$ is analytic when $\Re(\lambda)>\max(2\alpha_{0}, a)$. Thus, there is a real number
$\beta_0$ such that all roots of $H(\lambda)$ satisfy $\Re(\lambda) < \beta_0$ (refer to the discussion in
\cite[Lemma 4.1 in Chapter 1]{Hale93}), where
$$
\beta_0= \sup\{\Re(\lambda): H(\lambda) = 0, \lambda\in \C\}.
$$
Thus, for any $\beta>\beta_{0}$ there exists a positive constant
$C_{3}=C_{3}(\beta)$ such that
\begin{equation}\label{3.30}
\left|Y(t)\right|\leq C_{3}e^{\beta t},\quad t\geq 0.
\end{equation}

Now, we are ready to prove the conclusions.

1.  First, from Theorem \ref{Thm3.3} and (\ref{3.10}), there are two positive constants
$K_{3}$ and $K_{4}$ such that for $t\geq 0$,
\begin{eqnarray*}
P(t)&=&\left(\sigma_{0}+\sigma_{1}Ex(t;\phi)
+\sigma_{2}\int_{0}^{+\infty}K(s)Ex(t-s;\phi)ds \right)^{2}\\
&=&\left(\sigma_{0}+\sigma_{2}\int_{0}^{+\infty}K(s)Ex(t-s;\phi)ds \right)^{2}
+\sigma^{2}_{1}\left(Ex(t;\phi)\right)^{2}\\
&&{}+2\sigma_{1}Ex(t;\phi)\left(\sigma_{0}+\sigma_{2}\int_{0}^{+\infty}K(s)Ex(t-s;\phi)ds \right)\\
&\leq & \big(|\sigma_{0}|+|\sigma_{2}|\left(1+K_{1}\right) \|\phi\|\big)^{2}
+\sigma^{2}_{1}K^{2}_{1}\|\phi\|^{2}e^{2\alpha t}\\[0.1cm]
&&{}+2|\sigma_{1}| \Big(|\sigma_{0}|+|\sigma_{2}| \left(1+K_{1}\right) \|\phi\|\Big)K_{1}\|\phi\|e^{\alpha t}\\
&\leq & K_{3}+K_{4}e^{\alpha t},
\end{eqnarray*}
where
\begin{equation*}
K_{3}=\Big(|\sigma_{0}|+|\sigma_{2}|\left(1+K_{1}\right)\|\phi\|\Big)^{2}
\end{equation*}
and
\begin{equation*}
K_{4}=\sigma^{2}_{1}K^{2}_{1}\|\phi\|^{2}+2|\sigma_{1}| K_{1}\|\phi\| \left( |\sigma_{0}|+|\sigma_{2}|
\left(1+K_{1}\right)\|\phi\| \right).
\end{equation*}
We note that $K_3$ and $K_4$ are of order $\|\phi\|^2$.

If $\beta_{0} <0$, for any  $ \beta\in (\beta_{0}, 0)$  there exists a constant $C_{3}$ as
in (\ref{3.30}) such that for $t\geq 0,$
\begin{eqnarray*}
\left|M(t)\right|&=&\left|\int_{0}^{t}Y(s)P(t-s)ds\right|
\leq \int_{0}^{t}C_{3}e^{\beta s}\left( K_{3}+K_{4}e^{\alpha (t-s)}\right)ds\\
&=& C_{3}K_{3}\int_{0}^{t}e^{\beta s}ds+ C_{3}K_{4}\int_{0}^{t}e^{\beta s}e^{\alpha (t-s)}ds\\
&=&\dfrac{C_{3}K_{3}\left(1-e^{\beta t}\right)}{|\beta|}
+\dfrac{C_{3}K_{4}(e^{\beta t}-e^{\alpha t})}{\alpha-\beta}\\
&\leq & \dfrac{2C_{3}K_{3}}{|\beta|}
+\dfrac{2C_{3}K_{4}}{\left|\alpha-\beta\right|}.
\end{eqnarray*}
Therefore, the second moment is bounded for any initial function $\phi\in\left((-\infty,0],\R\right)$.

Now, let
\begin{equation*}
M_{\infty}=K_{3}\int_{0}^{+\infty}Y(t)dt,
\end{equation*}
then
\begin{eqnarray*}
\left|M(t)-M_{\infty}\right|
&=& \left|\int_{0}^{t}Y(s)P(t-s)ds-K_{3}\int_{0}^{+\infty}Y(s)ds \right|\\
&=&\left|\int_{0}^{t}Y(s)\left(P(t-s)-K_{3}\right)ds-K_{3}\int_{t}^{+\infty}Y(s)ds \right|\\
&\leq &K_{4}\int_{0}^{t}|Y(s)|e^{\alpha(t-s)}ds+K_{3}\int_{t}^{+\infty}|Y(t)|dt\\
&\leq & K_{3}\int_{t}^{+\infty}C_{3}e^{\beta t}dt+C_{3}K_{4}\int_{0}^{t}e^{\beta s}e^{\alpha(t-s)}ds\\
&\leq & \dfrac{C_{3}K_{3}e^{\beta t}}{|\beta|}+\dfrac{C_{3}K_{4}(e^{\beta t}-e^{\alpha t})}{\alpha-\beta}\\
&\leq &C_3\left(\dfrac{K_{3}}{|\beta|}+\dfrac{2K_{4}}{|\alpha-\beta|}\right)e^{t\max(\alpha , \beta)}.
\end{eqnarray*}
Thus, there exists a positive constant $C_{4}=C_3(\dfrac{K_{3}}{|\beta|}+\dfrac{2K_{4}}{|\alpha-\beta|})$ such that
\begin{eqnarray*}
\left|M(t)-M_{\infty}\right|
&\leq &C_{4}e^{t\max(\alpha , \beta)}\to 0 \quad \mathrm{as}\ t\to +\infty
\end{eqnarray*}
since $\max(\alpha , \beta)<0$, i.e., $M(t)$ approaches to $M_{\infty}$ exponentially
as $t\to +\infty$.

2. Assume $\beta_0 > 0$. We only need to show that there is a special solution $x(t;\phi)$ of \eqref{1.1} such that the corresponding second moment is unbounded. Similar to the proof of Theorem \ref{Thm3.6}, let $\lambda = \alpha + i \beta$ be a solution of $h(\lambda) = 0$, then $x_{\phi}(t) = \Re(e^{\lambda t})$ is a solution of \eqref{2.1} with initial function $\phi(\theta) = \Re(e^{\lambda \theta})\ (\theta \leq 0)$. Hence, for the solution $x(t;\phi)$  of \eqref{1.1} with this particular initial function, we have
$$P(t)= (\sigma_0 + e^{\alpha t} \Re[e^{i\beta t} (\sigma_1 + \sigma_2 \mathcal{L}(K)(\lambda))])^2.$$
Since the condition \textbf{H} is not satisfied, we have either $\sigma_0\not=0$ or $\sigma_1 + \sigma_2 \mathcal{L}(K)(\lambda)\not=0$. Thus, the function $P(t)\not\equiv 0$, and hence the Laplacian $\mathcal{L}(P)(s)$ is nonzero.

Since
$$M(t) = \mathcal{L}^{-1}\left(\frac{\mathcal{L}(P)(s)}{H(s)}\right) = \frac{1}{2\pi i} \lim_{T\to +\infty} \int_{c-i T}^{c+ i T} e^{s t} 
\frac{\mathcal{L}(P)(s)}{H(s)} d s$$
with $c > \beta_0$. Similar to the proof of Theorem \ref{Thm2.3}, and note that $\mathcal{L}(P)(s)$ is analytic when $\Re(s) = c>0$, there is
$\bar{\beta}\in (0,\beta_0)$ and a sequence $\{t_k\}$ with $t_k\to+\infty$ such that $M(t_k) > e^{\bar{\beta}t_k}$, which implies that the second 
moment is unbounded.
\qed

\begin{remark}\label{remark3.13}
The critical case when $\beta_{0}=0$ is not considered here, and the issue of boundedness criteria remains open.
\end{remark}

\section{Applications}

The functions $f_{1}(\lambda)$ and  $f_{2}(\lambda)$ in $H(\lambda)$ depend not only on the coefficients
of equation \eqref{1.1}, but also on the Laplace transforms of $X^{2}(t),\,X_{s}(t)X_{l}(t)$ and the delay kernel $K.$
Though it is possible to calculate these two functions numerically according to Lemma \ref{lem2.6},
it is not trivial to obtain $\beta_0 = \sup\{\Re(\lambda): H(\lambda) = 0\}$ for a given equation.
Hence the boundedness criteria established in Theorem \ref{Thm3.8} are not practical in applications. In applications, one need to
derive useful criteria according to the density kernel $K$ and coefficients of the equation \eqref{1.1}. Here, we give some
practical conditions, according to Theorem \ref{Thm3.8}, for applications.

First, in the case of discrete delay ($K(s)=\delta (s-1)$), we have
$$
f_{1}(\lambda) = \dfrac{\mathcal{L}(XX_{1})(\lambda)}{\mathcal{L}(X^2)(\lambda)} = g(\lambda, 1,0)
$$
and
$$
f_{2}(\lambda) = \dfrac{\mathcal{L}(X^2_{1})(\lambda)}{\mathcal{L}(X^2)(\lambda)}
=e^{-\lambda}.
$$
Thus
\begin{equation*}
H(\lambda)=\lambda-\left(2a+\sigma^{2}_{1}\right)-2\left(b+\sigma_{1}\sigma_{2}\right)g(\lambda,1,0)
-\sigma^{2}_{2}e^{-\lambda},
\end{equation*}
which give the same characteristic function $H(s)$ as in  \cite[Theorem 3.6]{LeiM} (here, we have provided explicit expressions for the functions $f(s)$ and $g(s)$ in \cite{LeiM} by \eqref{2.24}).
Therefore, sufficient conditions for the boundedness or unboundedness of the second moment
can be referred to \cite{LeiM}.

Next, if $b=0$, the fundamental solution $X(t)$ is known to us. Therefore, it is possible to obtain explicit sufficient conditions for the boundedness. Here we give a sufficient condition for the second moment to be bounded when $b=0$ and $K(s)=re^{-rs}$ ($j=1$ in the gamma distribution \eqref{2.5}).

\begin{theorem}\label{Thm4.1}
Let  $b=0$ and $K(s)=re^{-rs}\,(r>0).$  If $a<0$ and
\begin{equation}\label{4.1}
a_1>0,\quad a_3>0,\quad a_1a_2 - a_3 > 0
\end{equation}
where
\begin{eqnarray*}
a_1 &=& 3 (r-a) - \sigma_1^2,\\
a_2 &=& 2 r (r-a) - (2 a + \sigma_1^2) (3 r - a) - 2 r \sigma_1 \sigma_2,\\
a_3 &=& -2 a r(2 (r-a) - \sigma_1^2) - 2 r^2 (\sigma_1 + \sigma_2)^2,
\end{eqnarray*}
then the second moment is bounded.
\end{theorem}
\Proof If $b=0,$ the fundamental solution of \eqref{2.1} is given by
\begin{equation*}
X(t)=\left\{ \begin{array}{cc}
e^{at},\quad & t \geq 0,\\ 0,\quad & t<0
\end{array}\right.
\end{equation*}
and $\alpha_{0}=a$. Hence
the trivial solution of (\ref{2.1}) is locally asymptotically stable if and only if
$a<0$.

From the fundamental solution, we have for $\Re(\lambda)>2a,$
\begin{equation*}
\mathcal{L}(X^2)(\lambda)
=\int_0^\infty e^{-\lambda t} X^{2}(t)dt
=\int_0^\infty e^{-(\lambda-2a) t} dt
=\dfrac{1}{\lambda-2a},
\end{equation*}
\begin{equation*}
\mathcal{L}(XX_{s})(\lambda)
=\int_0^\infty e^{-\lambda t} X(t)X(t-s)dt
=e^{-(\lambda-a)s}\int_s^\infty e^{-(\lambda-2a) (t-s)} dt
=\dfrac{e^{-(\lambda-a)s}}{\lambda-2a},
\end{equation*}
and
\begin{eqnarray*}
\mathcal{L}(X_{s}X_{l})(\lambda)
&=&\int_0^\infty e^{-\lambda t} X(t-s)X(t-l)dt
=e^{-a(s+l)}\int_{s\vee l}^\infty e^{-(\lambda-2a)t} dt\\
&=&\left\{ \begin{array}{cc}
\dfrac{e^{-(\lambda-a)s}e^{-al}}{\lambda-2a},\quad & s \geq l\geq 0,\\[0.3cm]
\dfrac{e^{-(\lambda-a)l}e^{-as}}{\lambda-2a},\quad & 0 \leq s < l,
\end{array}\right.
\end{eqnarray*}
where $s\vee l=\max\{s, l\}$. Therefore
$$
g(\lambda, s,l)=\dfrac{\mathcal{L}(X_{s}X_{l})(\lambda)}{\mathcal{L}(X^2)(\lambda)}
=\left\{ \begin{array}{cc}
e^{-(\lambda-a)s}e^{-al},\quad & s \geq l\geq 0,\\[0.1cm]
e^{-(\lambda-a)l}e^{-as},\quad & 0 \leq s < l.
\end{array}\right.
$$

Let $K(s)=re^{-rs}\,(r>0)$, then
$$
f_{1}(\lambda) = \dfrac{r}{\lambda+r-a},\quad f_2(\lambda)=\dfrac{2r^2}{(\lambda+r-a)(\lambda+2r)}.
$$
Thus from (\ref{3.7}),
\begin{equation}\label{4.2}
H(\lambda)=\lambda-\left(2a+\sigma^{2}_{1}\right)-\dfrac{2r\sigma_{1}\sigma_{2}}{\lambda+r-a}
-\dfrac{2r^2\sigma^{2}_{2}}{(\lambda+r-a)(\lambda+2r)}.
\end{equation}
Hence
$H(\lambda)=0$ if and only if
\begin{equation}\label{4.3}
\bar{H}(\lambda)=\lambda^3+a_{1}\lambda^2+a_{2}\lambda+a_{3}=0,
\end{equation}
where
\begin{eqnarray*}
a_{1}&=&3(r-a)-\sigma^2_{1},\; a_{2}=2r(r-a)-(2a+\sigma^2_{1})(3r-a)-2r\sigma_{1}\sigma_{2},\\
a_{3}&=&-2ar(2(r-a)-\sigma^2_{1})-2r^2(\sigma_{1}+\sigma_{2})^2.
\end{eqnarray*}
From the Routh-Hurwitz criterion, all roots of $\bar{H}(\lambda)=0$ have negative real
parts if and only if
$$a_{1}>0,\; a_{3}>0 \;\; \text{and} \;\;a_{1}a_{2}-a_{3}>0.$$
Thus, the theorem is proved.
\qed

The following theorem gives a sufficient condition for the unboundedness of the second moment for general situations.
\begin{theorem}\label{Thm4.2}
If either
\begin{equation}\label{4.4}
b+\sigma_{1}\sigma_{2}\leq 0,\quad \sigma^{2}_{2}-2\sigma_{1}\sigma_{2}-\sigma^{2}_{1}<2(a+b)
\end{equation}
or
\begin{equation}\label{4.5}
b+\sigma_{1}\sigma_{2}\geq 0,\quad \sigma^{2}_{2}+2\sigma_{1}\sigma_{2}-\sigma^{2}_{1}<2(a-b),
\end{equation}
the second moment is unbounded.
\end{theorem}
\Proof From \eqref{2.22}, we have
$$|g(0,s,l)| \leq 1,$$
and therefore
$$
|f_{1}(0)|\leq \int_{0}^{+\infty}K(s)ds \leq 1,\quad
|f_{2}(\lambda)|\leq (\int_{0}^{+\infty}K(s)ds)^2 \leq 1.
$$
Thus, we have
$$
H(0)= -(2a+\sigma^{2}_{1}) -2(b+\sigma_{1}\sigma_{2})f_{1}(0)-\sigma^{2}_{2}f_{2}(0) < 0
$$
when either \eqref{4.4} or \eqref{4.5} is satisfied.
Furthermore, it is easy to have $H(\lambda)>0$ when $\lambda\in \R$ is large enough. Thus the equation
$H(\lambda)=0\,(\lambda \in \R)$ has at least one positive solution, which implies that
$\beta_{0}>0$, and the second moment is unbounded by Theorem \ref{Thm3.8}.
\qed

\section{An example}
\begin{figure}[htbp]
\centering
\includegraphics[width=7cm]{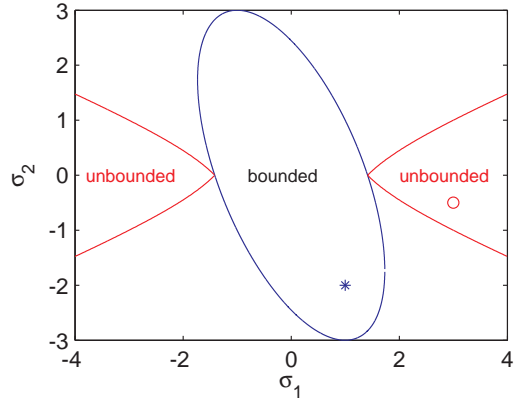}
\caption{The bounded and unbounded regions of the second moment obtained from Theorems \ref{Thm4.1} and \ref{Thm4.2},
where $r=0.5.$ }
\label{fig1}
\end{figure}

\begin{figure}[htbp]
\centering
\includegraphics[width=12cm]{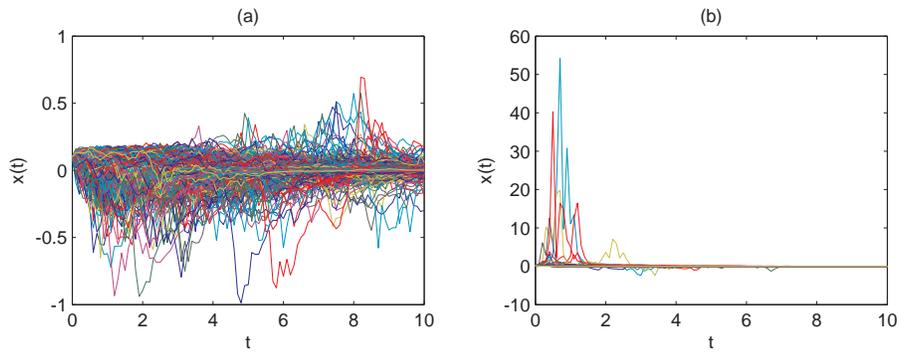}
\caption{Numerical results of 1000 sample solutions of \eqref{4.1}. Parameters used are (a) $(\sigma_1, \sigma_2) = (1,-2)$ (the star in Fig. \ref{fig1}), and (b) $(\sigma_1,\sigma_2) = (3, -0.5)$ (the circle in Fig. \ref{fig1}). All initial functions
are taken as $x(t) = 0.1$ for $t < 0$.}
\label{fig2}
\end{figure}

Here, we consider an example of following linear stochastic delay differential equation
\begin{equation}\label{5.1}
dx(t)= -x(t)dt+\left(\sigma_{1}x(t)+\sigma_{2}\int_{0}^{+\infty}K(s)x(t-s)ds\right)dW_t,
\end{equation}
where $K(s)=re^{-rs}\ (r>0)$. Fig. \ref{fig1} shows regions in the $(\sigma_1,\sigma_2)$ plane to have bounded and unbounded second moments according to Theorem \ref{Thm4.1} and Theorem \ref{Thm4.2}, respectively. Fig. \ref{fig2} show sample solutions with $(\sigma_1, \sigma_2) = (1,-2)$ (the star in Fig. \ref{fig1}) and with $(\sigma_1,\sigma_2) = (3, -0.5)$ (the circle in Fig. \ref{fig1}), respectively. Simulations show that when $(\sigma_1, \sigma_2) = (1,-2)$, all sample solutions are bounded. But when $(\sigma_1,\sigma_2) = (3, -0.5)$, the sample solutions have a positive probability to reach a large value. These numerical results show agreement with our theoretical analysis.

\vspace{0.5cm}

\noindent\textbf{Final Remark.} \textit{All results in this paper are obtained under the It\^{o} interpretation.
Analogous results can be obtained for the Stratonovich interpretation.}

\begin{appendix}
\section{Proofs of Lemmas \ref{lem2.1} and \ref{lem2.2}}

\noindent\textit{Proof of Lemma \ref{lem2.1}.} Since $x_{\phi}(t)$ satisfies (\ref{2.1}) for $t\geq 0$,
\begin{equation*}
x_{\phi}(t)= \phi(0)+ a\int_{0}^{t}x_{\phi}(s)ds+ b\int_{0}^{t}
\int_{0}^{+\infty}K(\theta)x_{\phi}(s-\theta)d\theta ds, \quad t\geq 0
\end{equation*} and $x_{\phi}(t)=\phi(t)$ for $t\in (-\infty, 0].$ Therefore for $t\geq 0$
\begin{eqnarray*}
|x_{\phi}(t)| &\leq & |\phi(0)|+ |a|\int_{0}^{t}|x_{\phi}(s)|ds
+|b|\int_{0}^{t}\int_{s}^{+\infty}K(\theta)|\phi(s-\theta)|d\theta ds\\
&&{}+|b|\int_{0}^{t}\int_{0}^{s}K(\theta)|x_{\phi}(s-\theta)|d\theta ds\\
&\leq &\|\phi\|+ |a|\int_{0}^{t}|x_{\phi}(s)|ds
+|b|\cdot\|\phi\|\int_{0}^{t}\int_{s}^{+\infty}K(\theta)d\theta ds\\
&&{}+|b|\int_{0}^{t}|x_{\phi}(\theta)|\int_{0}^{t-\theta}K(v)dv d\theta .
\end{eqnarray*}
Then from (\ref{2.2}) and (\ref{2.3}), we have
\begin{eqnarray*}
|x_{\phi}(t)|
&\leq &\|\phi\|+ |a|\int_{0}^{t}|x_{\phi}(s)|ds\\
&&{}
+|b|\cdot\|\phi\|\int_{0}^{t}\int_{s}^{+\infty}K(\theta)d\theta ds
+|b|\int_{0}^{t}|x_{\phi}(\theta)|d\theta\\
&\leq &\|\phi\|+ \left(|a|+|b|\right)\int_{0}^{t}|x_{\phi}(s)|ds
+|b|\cdot\|\phi\|\int_{0}^{t}e^{-\mu s}\int_{s}^{+\infty}e^{\mu \theta}K(\theta)d\theta ds\\
&\leq &\|\phi\|+ \left(|a|+|b|\right)\int_{0}^{t}|x_{\phi}(s)|ds
+\rho|b|\cdot\|\phi\|\int_{0}^{t}e^{-\mu s} ds\\
&\leq &\|\phi\|\left(1+\dfrac{\rho|b|}{\mu}\right)+ \left(|a|+|b|\right)
\int_{0}^{t}|x_{\phi}(s)|ds.
\end{eqnarray*}
Using the Gronwall inequality, we obtain
\begin{equation*}
|x_{\phi}(t)|\leq \|\phi\|\left(1+\dfrac{\rho|b|}{\mu}\right)e^{\left(|a|+|b|\right)t}, \quad t\geq 0.
\end{equation*}
Thus the inequality (\ref{2.2}) is satisfied with
$A=\|\phi\|\left(1+\dfrac{\rho|b|}{\mu}\right)$ and  $ B=\left(|a|+|b|\right).$
\qed

\vspace{0.5cm}

\noindent\textit{Proof of Lemma \ref{lem2.2}.}
 Taking Laplace transform to both sides of (\ref{2.1}), we obtain
\begin{equation*}
\int_{0}^{+\infty}e^{-\lambda t}\dfrac{dx_{\phi}(t)}{dt}dt
= a\int_{0}^{+\infty}e^{-\lambda t}x_{\phi}(t)dt
+ b\int_{0}^{+\infty}e^{-\lambda t}\int_{0}^{+\infty}K(s)x_{\phi}(t-s)ds dt,
\end{equation*}
when $t\geq 0$, i.e.,
\begin{eqnarray*}
-\phi(0)+\lambda\mathcal{L}(x_{\phi})&=& a\mathcal{L}(x_{\phi})
+ b\int_{0}^{+\infty}e^{-\lambda t}\int_{t}^{+\infty}K(s)\phi(t-s)ds dt\\
&&+ b\int_{0}^{+\infty}e^{-\lambda t}\int_{0}^{t}K(s)x_{\phi}(t-s)ds dt, \,\, (t\geq 0).
\end{eqnarray*}
Since
\begin{eqnarray*}
\int_{0}^{+\infty}e^{-\lambda t}\int_{0}^{t}K(s)x_{\phi}(t-s)ds dt
&=&\int_{0}^{+\infty}K(s)\int_{s}^{+\infty}e^{-\lambda t}x_{\phi}(t-s)dt ds\\
&=&\int_{0}^{+\infty}K(s)e^{-\lambda s}\int_{0}^{+\infty}e^{-\lambda l}x_{\phi}(l)dl ds\\
&=&\mathcal{L}(K)\mathcal{L}(x_{\phi}),
\end{eqnarray*}
we have
\begin{equation*}
\left(\lambda-a-b\mathcal{L}(K)\right)\mathcal{L}(x_{\phi})=
\phi(0)+b\int_{0}^{+\infty}e^{-\lambda t}\int_{t}^{+\infty}K(s)\phi(t-s)ds dt,
\end{equation*}
which yields,
\begin{equation*}
\mathcal{L}(x_{\phi})=
\dfrac{\phi(0)}{h(\lambda)}+\dfrac{b}{h(\lambda)}
\int_{0}^{+\infty}e^{-\lambda t}\int_{t}^{+\infty}K(s)\phi(t-s)ds dt.
\end{equation*}
Therefore, since $\mathcal{L}(X)=h^{-1}(\lambda)$, we have
\begin{eqnarray*}
x_{\phi}(t)&=&
\mathcal{L}^{-1}\left(\dfrac{1}{h(\lambda)}\right)\phi(0)
+b\mathcal{L}^{-1}\left(\dfrac{1}{h(\lambda)}
\int_{0}^{+\infty}e^{-\lambda t}\int_{t}^{+\infty}K(s)\phi(t-s)ds dt\right)\\
&=&X(t)\phi(0)+b\int_{0}^{t}X(t-s)\int_{s}^{+\infty}K(\theta)\phi(s-\theta)d\theta ds,
\end{eqnarray*}
and the Lemma is proved.
\qed

\vspace{0.5cm}

\section{Proof of Lemma \ref{lem2.7}}

First, when $\Re(\lambda)>\max\{2\alpha_0,a,a+\alpha_0\}$, since
\begin{eqnarray}
&&\frac{1}{2\pi}\int_{-\infty}^\infty \frac{1}{h(i\omega)h(\lambda - i\omega -a)}d \omega \nonumber\\
&=& \frac{1}{2\pi} \int_{-\infty}^\infty \frac{1}{h(i\omega)}\int_0^{+\infty}e^{-(\lambda - i\omega - a)t}d t d\omega\nonumber\\
&=&\int_0^{+\infty}e^{-(\lambda - a)t}\frac{1}{2\pi}\int_{-\infty}^\infty \frac{e^{i\omega t}}{h(i\omega)}d\omega dt\nonumber\\
\label{eq:B0}
&=&\int_0^{+\infty} e^{-(\lambda - a)t} X(t) d t = \frac{1}{h(\lambda - a)},
\end{eqnarray}
then
\begin{eqnarray}\label{B.1}
\mathcal{L}(X^2)(\lambda)
&=&\dfrac{1}{2\pi}\int_{-\infty}^\infty \dfrac{1}{h(i\omega) (\lambda- i \omega - a)} d \omega \nonumber\\
&&{}
+ \dfrac{1}{2\pi} \int_{-\infty}^\infty \dfrac{1}{h(i\omega)}\left(\dfrac{1}{h(\lambda-i\omega)}
- \dfrac{1}{\lambda - i\omega - a}\right) d\omega \nonumber \\
&=&\dfrac{1}{h(\lambda-a)} + \dfrac{1}{2\pi}\int_{-\infty}^{\infty} \dfrac{b\mathcal{L}(K)(\lambda-i\omega)}
{h\left(i\omega\right)h\left(\lambda-i\omega\right) \left(\lambda - i\omega - a\right)}d \omega \nonumber\\
&=& \dfrac{1}{h(\lambda-a)}\left(1+g(\lambda)\right),
\end{eqnarray}
where
\begin{equation}\label{B.2}
g(\lambda)=\dfrac{h(\lambda-a)}{2\pi}\int_{-\infty}^{\infty} \dfrac{b\mathcal{L}(K)(\lambda-i\omega)}
{h\left(i\omega\right)h\left(\lambda-i\omega\right) \left(\lambda - i\omega - a\right)}d \omega
\end{equation}
is convergent for $\Re(\lambda)>\max\{2\alpha_0,a,a+\alpha_0\}$.  We have the following result for $g(\lambda)$.

\noindent\textbf{Lemma B.1.} \textit{
Let $g(\lambda)$ be defined as in \eqref{B.2}. Then for any $\Re(\lambda) > \max\{2\alpha_{0}, a,a+\alpha_0,-\mu\}$,
\begin{equation*}
\lim_{|\lambda|\rightarrow +\infty}\left|g(\lambda)\right|=0.
\end{equation*}
}

\Proof First, when $\Re(\lambda)>-\mu$,
\begin{eqnarray}\label{B.3}
\left|\mathcal{L}(K)(\lambda-i\omega)\right|
&=&\left|\int_{0}^{+\infty}e^{-(\lambda-i\omega)t}K(t)dt\right|
\leq \int_{0}^{+\infty}e^{-Re(\lambda)t}K(t)dt \nonumber\\
&\leq & \int_{0}^{+\infty}e^{\mu t}K(t)dt=\rho,
\end{eqnarray}
and hence,
$$
|g(\lambda)|\leq \frac{|b|\rho }{2 \pi} \int_{-\infty}^\infty\frac{|h(\lambda - a)|}{|h(i\omega)h(\lambda-i\omega)(\lambda-i\omega-a)|}d\omega.
$$
Given a positive constant $\omega_{0}$ such that
$\omega_{0}> |\lambda| +|a|+|b|\;(\lambda\in\C)$. Then for any $|\omega|> \omega_{0},$
\begin{equation*}
\dfrac{1}{\left|h(\lambda - i\omega)\right|}
\leq \dfrac{1}{|\omega |-|\lambda|-|a | - |b|}.
\end{equation*}
Thus for any $|\omega|> \omega_{0},$
\begin{eqnarray*}
\left|\dfrac{1}{h(i\omega)h\left(\lambda-i\omega\right)\left(\lambda - i\omega - a \right)}\right|
&\leq & \dfrac{1}{\left(\left|\omega \right| - \left|a \right| - \left|b \right|\right)
\left(\left|\omega \right| -\left|\lambda\right|- \left|a \right| - \left|b \right|\right)
\left(\left|\omega \right| -\left|\lambda\right|- \left|a \right|\right) }\\
&\leq &\dfrac{1}{\left(\left|\omega \right|-\left|\lambda\right|
- \left|a \right| - \left|b \right|\right)^{3}}.
\end{eqnarray*}
Therefore when $\Re(\lambda)> -\mu$,
\begin{eqnarray*}
|g(\lambda)|&\leq& \frac{|b|\rho}{2\pi}\int_{-\infty}^{-\omega_0} \frac{|h(\lambda-a)|}{|h(i\omega)h(\lambda-i\omega)(\lambda-i\omega-a)|}d\omega\\
&&{} + \frac{|b|\rho}{2\pi}\int_{-\omega_0}^{\omega_0} \frac{|h(\lambda-a)|}{|h(i\omega)h(\lambda-i\omega)(\lambda-i\omega-a)|}d\omega\\
&&{} + \frac{|b|\rho}{2\pi}\int_{\omega_0}^{\infty} \frac{|h(\lambda-a)|}{|h(i\omega)h(\lambda-i\omega)(\lambda-i\omega-a)|}d\omega\\
&\leq& \frac{|b|\rho}{2\pi}\times 2\int_{\omega_0}^\infty \frac{|h(\lambda - a)|}{(\omega - |\lambda| - |a| - |b|)^3} d \omega\\
&&{} +\frac{|b|\rho}{2\pi}\int_{-\omega_0}^{\omega_0} \frac{|h(\lambda-a)|}{|h(i\omega)h(\lambda-i\omega)(\lambda-i\omega-a)|}d\omega\\
&=& \frac{|b|\rho }{2\pi}\Big(\frac{|h(\lambda -a)|}{(\omega_0 - |\lambda| - |a| - |b|)^2}\\
&&{}\qquad\quad+\int_{-\omega_0}^{\omega_0} \frac{|h(\lambda-a)|}{|h(i\omega)h(\lambda-i\omega)(\lambda-i\omega-a)|}d\omega\Big)
\end{eqnarray*}
Now, since
$$0\leq \lim_{|\lambda|\to +\infty}\frac{|h(\lambda - a)|}{(\omega_0 - |\lambda| - |a| - |b|)^2} \leq \lim_{|\lambda| + \infty} \frac{|\lambda| + 2 |a| + |b| \rho}{(\omega_0 - |\lambda| - |a| - |b|)^2} = 0,$$
and when $\Re(\lambda) > \max\{2 \alpha_0, a,-\mu\}$,
\begin{eqnarray*}
0&\leq& \lim_{|\lambda|\to+\infty}\int_{-\omega_0}^{\omega_0} \frac{|h(\lambda-a)|}{|h(i\omega)h(\lambda-i\omega)(\lambda-i\omega-a)|}d\omega \\
&=& 2 \int_0^{\omega_0}\lim_{|\lambda|\to+\infty}\frac{|h(\lambda-a)|}{|h(i\omega)|h(\lambda-i\omega)||(\lambda-i\omega-a)|} = 0,
\end{eqnarray*}
we have
$$\lim_{|\lambda|\to+\infty}|g(\lambda)| = 0,\quad \Re(\lambda) > \max\{2\alpha_0, a, a+\alpha_0,-\mu\}.$$
The lemma is proved.
\qed

\noindent\textit{Proof of Lemma \ref{lem2.7}.} Similar to \eqref{eq:B0}, we have, when $s\geq l$,
$$\frac{e^{-\lambda s}}{2\pi} \int_{-\infty}^\infty \frac{e^{i\omega (s-l)}}{h(i\omega)(\lambda-i\omega-a)}d\omega = \frac{e^{-(\lambda-a)}e^{-a s}}{h(\lambda - a)}.$$
Thus, from \eqref{2.23},
when $\Re(\lambda)>\max\{2\alpha_0,a\}$, we obtain
\begin{eqnarray}\label{B.7}
\mathcal{L}(X_{s}X_{l})
&=&\dfrac{e^{-\lambda s}}{2\pi}\int_{-\infty}^\infty \dfrac{e^{i\omega (s-l)}}{h(i\omega) (\lambda- i \omega-a )} d \omega \nonumber\\
&&{}
+ \dfrac{e^{-\lambda s}}{2\pi} \int_{-\infty}^\infty \dfrac{e^{i\omega (s-l)}}{h(i\omega)}\left(\dfrac{1}{h(\lambda-i\omega)}
- \dfrac{1}{\lambda - i\omega-a }\right) d\omega \nonumber\\
&=&\dfrac{e^{-(\lambda-a) l}e^{-as}}{h(\lambda-a)}
+ \dfrac{e^{-\lambda s}}{2\pi}\int_{-\infty}^\infty\dfrac{be^{i\omega (s-l)}\mathcal{L}(K)(\lambda-i\omega)}
{h(i\omega)h(\lambda-i\omega)\left(\lambda - i\omega-a\right)} d \omega \nonumber\\
&=& \dfrac{1}{h(\lambda-a)}\left(e^{-(\lambda-a) l}e^{-as}+e^{-\lambda s}\tilde{g}_1(\lambda,s,l)\right),
\end{eqnarray}
where
\begin{equation}\label{B.8}
\tilde{g}_1(\lambda,s,l)= \dfrac{h(\lambda-a)}{2\pi}\int_{-\infty}^\infty\dfrac{be^{i\omega (s-l)}
\mathcal{L}(K)(\lambda-i\omega)}
{h(i\omega)h(\lambda-i\omega)\left(\lambda - i\omega-a\right)} d \omega.
\end{equation}
Similarly, when $0\leq s<l$,
\begin{eqnarray*}
\mathcal{L}(X_{s}X_{l})
&=& \dfrac{1}{h(\lambda-a)}\left(e^{-(\lambda-a) s}e^{-al}+e^{-\lambda l}\tilde{g}_{2}(\lambda,s,l)\right),
\end{eqnarray*}
where
\begin{equation}\label{B.9}
\tilde{g}_{2}(\lambda,s,l)= \dfrac{h(\lambda-a)}{2\pi}\int_{-\infty}^\infty\dfrac{be^{i\omega (l-s)}
\mathcal{L}(K)(\lambda-i\omega)}
{h(i\omega)h(\lambda-i\omega)\left(\lambda - i\omega-a\right)} d \omega.
\end{equation}
From (\ref{B.8}), (\ref{B.9}) and \eqref{B.3}, we obtain
\begin{equation*}
\left|\tilde{g}_{j}(\lambda,s,l)\right|\leq
\dfrac{\left|h(\lambda-a)\right|}{2\pi}\int_{-\infty}^\infty\dfrac{\left|b\right|\rho}
{\left|h(i\omega)h(\lambda-i\omega)\left(\lambda - i\omega\right)\right|} d \omega \,(j=1, 2).
\end{equation*}
Thus similar to the proof of Lemma B.1, for $\Re(\lambda)>\max\{2\alpha_0,a,-\mu\}$, we have
\begin{equation}
\lim_{|\lambda|\to +\infty}\left|\tilde{g}_{j}(\lambda,s,l)\right|=0 \;\;
(\text{uniformly with respect to $s$ and $l,$ }\,j=1,2).
\label{B.10}
\end{equation}
From \eqref{2.21} and \eqref{B.1}, we have for $\Re(\lambda) > \max\{ 2 \alpha_0,a, a+\alpha_0, -\mu\}$,
\begin{eqnarray}
g(\lambda, s,l) &=& \frac{\mathcal{L}(X_sX_l)(\lambda)}{\mathcal(X^2)(\lambda)}\nonumber\\
\label{eq:B11}
&=&\left\{
\begin{array}{ll}
\dfrac{e^{-(\lambda-a)l}e^{-as} + e^{-\lambda s} \tilde{g}_1(\lambda, s, l)}{1 + g(\lambda)},\quad & (s\geq l\geq 0),\\
\dfrac{e^{-(\lambda-a)s}e^{-as} + e^{-\lambda l} \tilde{g}_2(\lambda, s, l)}{1 + g(\lambda)}\quad & (0\leq s\leq l).
\end{array}
\right.
\end{eqnarray}
From Lemma B.1 and \eqref{B.10}, when $\Re(\lambda)>\max\{2\alpha_0,a,a+\alpha_0,-\mu\}$,
for any $\varepsilon>0$, there exists a constant $T_{0}=T_{0}(\varepsilon)$,
independent of $s$ and $l$ such that for $|\lambda|>T_{0}$,
\begin{equation*}
\left|g(\lambda)\right|< \varepsilon,\quad \mathrm{and}\quad
\left|\tilde{g}_{j}(\lambda,s,l)\right|<\varepsilon\quad (j=1,2).
\end{equation*}
Therefore for the above $\varepsilon$ and $T_{0}$, and $|\lambda|>T_{0}$, from \eqref{eq:B11}, we have
$$
\left|g(\lambda,s,l)\right|
\leq \left\{
\begin{array}{ll}
\dfrac{e^{-(T_{0}-a)l}e^{-as}}{1-\varepsilon}+\dfrac{\varepsilon e^{-T_{0}s}}{1-\varepsilon}, \quad & s\geq l\geq 0,\\
\dfrac{e^{-(T_0 - a)s}e^{-al}}{1-\varepsilon}+\dfrac{\varepsilon e^{-T_0l}}{1-\varepsilon}, \quad &0\leq s \leq l,
\end{array}\right.
$$
and hence
$$\lim_{|\lambda|\to+\infty} |g(\lambda,s, l)| = 0.$$
Thus, (\ref{2.25}) and  (\ref{2.26}) are satisfied and the Lemma is proved.
\qed

\end{appendix}

\section*{References}
\bibliographystyle{elsarticle-num}

\end{document}